\documentclass[12pt]{article}%
\usepackage{graphicx}
\usepackage{amsmath}
\usepackage{a4}
\usepackage{amsfonts}
\usepackage{amssymb}%
\providecommand{\U}[1]{\protect \rule{.1in}{.1in}}
\newtheorem{theorem}{Theorem}

\newtheorem{corollary}[theorem]{Corollary}

\newtheorem{lemma}[theorem]{Lemma}

\newtheorem{proposition}[theorem]{Proposition}
\newtheorem{remark}[theorem]{Remark}

\newenvironment{proof}[1][Proof]{\textbf{#1.} }{\  \rule{0.5em}{0.5em}}

\begin{document}

\title{On the Mapping class group of a genus $2$ handlebody.}
\author{Charalampos Charitos$\dagger$, Ioannis Papadoperakis$\dagger$\\and Georgios Tsapogas$\ddagger$\\$\dagger$Agricultural University of Athens \\and $\ddagger$University of the Aegean}
\maketitle

\begin{abstract}
A complex of incompressible surfaces in a handlebody is constructed so that it
contains, as a subcomplex, the complex of curves of the boundary of the
handlebody. For genus $2$ handlebodies, the group of automorphisms of this
complex is used to characterize the mapping class group of the handlebody. In
particular, it is shown that all automorphisms of the complex of
incompressible surfaces are geometric, that is, induced by a homeomorphism of
the handlebody.\newline2000 \textit{{Mathematics Subject Classification:}
57N10, 57N35 }

\end{abstract}

\section{Introduction and statements of results}

Let $M$ be a $3-$dimensional orientable differentiable manifold with or
without boundary and let $S\subset M$ be a properly embedded surface i.e., the
interior $\operatorname*{Int}\left(  S\right)  $ and the boundary $\partial S$
of $S$ satisfy the inclusions
\[
\operatorname*{Int}\left(  S\right)  \subset \operatorname*{Int}\left(
M\right)  \text{\  \ and\  \ }\partial S\subset \partial M,
\]
$S$ is transverse to $\partial M,$ and the intersection of $S$ with a compact
subset of $M$ is compact in $S.$ A \emph{compressible disk} for $S$ is an
embedded disk $D\subset M$ such that $\partial D\subset S,$
$\operatorname*{Int}\left(  D\right)  \subset M\setminus S$ and $\partial D$
is an essential loop in $S,$ i.e., the map $\partial D\rightarrow S$ induces
an injection $\pi_{1}\left(  \partial D\right)  \rightarrow \pi_{1}\left(
S\right)  .$ A properly embedded surface $S\subset M$ is \emph{incompressible}
if there are no compressible disks for $S$ and no component of $S$ is a sphere
that bounds a ball. Recall also that a map $F:X\times \left[  0,1\right]
\rightarrow Y$ is a proper isotopy if for all $t\in I,$ $F\bigm \vert_{X\times
\left \{  t\right \}  }$ is a proper embedding. In this case we will be saying
that $F\left(  X\times \left \{  0\right \}  \right)  $ and $F\left(
X\times \left \{  1\right \}  \right)  $ are properly isotopic in $Y.$

Let $M$ denote a handlebody of genus $n.$ The complex $\mathcal{C}\left(
\partial M\right)  $ of curves for (any surface and, hence, for) the boundary
$\partial M$ is defined (see Harvey \cite{Har}) as follows: the vertices are
isotopy classes of essential unoriented (in general, non-boundary parallel,
but this is irrelevant as $\partial M$ is a closed surface) simple closed
curves in $\partial M$ and the simplices of $\mathcal{C}\left(  \partial
M\right)  $ are $\left(  \left[  \alpha_{1}\right]  ,\ldots,\left[  \alpha
_{k}\right]  \right)  $ where $\left[  \alpha_{i}\right]  ,\left[  \alpha
_{j}\right]  $ are distinct classes having disjoint representatives for $i\neq
j.$ We similarly define the complex $\mathcal{I}\left(  M\right)  $ of
(properly embedded and connected) incompressible surfaces in $M$: a vertex
$\left[  S\right]  $ in $\mathcal{I}\left(  M\right)  $ is an isotopy class of
connected properly embedded surfaces in $M$ with the additional requirement
that, unless $S$ is a meridian, it is isotopic to a surface $\overline{S}$
embedded in $\partial M$ via an isotopy
\[
F:S\times \left[  0,1\right]  \rightarrow M
\]
with $F\left(  S\times \left \{  0\right \}  \right)  =S$ and $F\left(
S\times \left \{  1\right \}  \right)  =\overline{S}.$ Note that, as
$\overline{S}\subset \partial M$ is not properly embedded in $M,$ such an
isotopy cannot be proper. However, we require that $F$ is proper when
restricted to $\left[  0,1\right)  .$ The simplices of $\mathcal{I}\left(
M\right)  $ are $\left(  \left[  S_{1}\right]  ,\ldots,\left[  S_{k}\right]
\right)  $ where $\left[  S_{i}\right]  ,\left[  S_{j}\right]  $ are distinct
isotopy classes of surfaces having disjoint representatives for $i\neq j.$

Observe that not all incompressible surfaces properly embedded in $M$ are
isotopic to a surface embedded in $\partial M$ (see for example the high genus
surfaces constructed in \cite{Qiu},\cite{CPT}). In the sequel, by saying that
a properly embedded surface is isotopic to a surface contained in the boundary
we will mean that such an isotopy is proper when restricted to time interval
$\left[  0,1\right)  .$

We consider two special classes of incompressible surfaces in $M,$ namely,
meridians in $M$ and properly embedded surfaces in $M$ which are homeomorphic
to an annulus with the two boundary components being isotopic on $\partial M.$
These two classes of surfaces define two subcomplexes of $\mathcal{I}\left(
M\right)  $ in a similar way: the vertices are isotopy classes of meridians
(resp. annuli) and the simplices collections of pair-wise distinct classes of
meridians (resp. annuli) having disjoint representatives. We call the first
one the (well known) complex of meridians $\mathcal{D}\left(  M\right)  $ and
the second the complex of annuli $\mathcal{A}\left(  M\right)  .$ Every
essential simple closed curve on $\partial M$ which bounds a disk in $M$ gives
rise to a meridian which is trivially an incompressible surface. Moreover, two
meridians are disjoint up to isotopy if and only if their boundaries are
disjoint up to isotopy. Therefore, the subcomplex of meridians $\mathcal{D}%
\left(  M\right)  \subset \mathcal{I}\left(  M\right)  $ can be viewed as a
subcomplex of the complex of curves $\mathcal{C}\left(  \partial M\right)  .$

The complex of meridians $\mathcal{D}\left(  M\right)  $ has been introduced
in \cite{MC} where it was used in the study of the mapping class groups of
$3-$ manifolds. It has also been studied in \cite{M-M} where it is shown to be
a quasi-convex subset of $\mathcal{C}\left(  \partial M\right)  ,$ $M$ being a
$3-$manifold with boundary.

Similarly, if $\alpha$ is an essential simple closed curve on $\partial M$
which is not homotopically trivial in $M,$ we may consider a properly embedded
annular surface $S_{\alpha}$ whose boundary consists of two (parallel) copies
of $\alpha$ so that $S_{\alpha}$ is isotopic to the annular region in
$\partial M$ bounded by its boundary components. Again, two such surfaces
$S_{\alpha},S_{\beta}$ are distinct up to isotopy if and only if the
corresponding curves $\alpha,\beta$ are distinct, up to isotopy, on $\partial
M.$ Therefore, the subcomplex of annuli $\mathcal{A}\left(  M\right)
\subset \mathcal{I}\left(  M\right)  $ can be viewed as a subcomplex of the
complex of curves $\mathcal{C}\left(  \partial M\right)  .$ Moreover, we have
a bijection
\[
\mathcal{C}\left(  \partial M\right)  \longleftrightarrow \mathcal{D}\left(
M\right)  \cup \mathcal{A}\left(  M\right)  .
\]
We will be writing $\mathcal{C}\left(  \partial M\right)  $ to denote the
subcomplex $\mathcal{D}\left(  M\right)  \cup \mathcal{A}\left(  M\right)  $ in
$\mathcal{I}\left(  M\right)  .$

Our goal is to show that when $M$ is a handlebody of genus $2$ the
automorphisms of the complex of incompressible surfaces are all geometric,
that is, they are induced by homeomorphisms of $M.$ However, this is not true
for the whole complex $\mathcal{I}\left(  M\right)  .$ A non-geometric
automorphisms of the complex $\mathcal{I}\left(  M\right)  $ can be seen as
follows: there exist non-separating curves $\alpha,$ $\beta,$ $\gamma$ in
$\partial M$ such that the closures of the components of $\partial
M\setminus \left(  \alpha \cup \beta \cup \gamma \right)  $ are pairs of pants
$P_{\alpha,\beta,\gamma}^{+},$ $P_{\alpha,\beta,\gamma}^{-}$ which are non
isotopic (see Remark \ref{r5}). The corresponding vertices $\left[
P_{\alpha,\beta,\gamma}^{+}\right]  ,$ $\left[  P_{\alpha,\beta,\gamma}%
^{-}\right]  $ in $\mathcal{I}\left(  M\right)  $ are of finite valence and
there is an infinite number of (pairs of) vertices of this type in
$\mathcal{I}\left(  M\right)  .$ It is shown (cf. proof of Theorem
\ref{main_theorem}) that the automorphism of $\mathcal{I}\left(  M\right)  $
which interchanges $\left[  P_{\alpha,\beta,\gamma}^{+}\right]  ,$ $\left[
P_{\alpha,\beta,\gamma}^{-}\right]  $ and fixes all other vertices cannot be geometric.

Consider the subcomplex $\mathcal{I}_{0}\left(  M\right)  $ of $\mathcal{I}%
\left(  M\right)  $ consisting of all vertices of infinite valence and we show
that the map
\[
A_{0}:\mathcal{MCG}\left(  M\right)  \rightarrow Aut\left(  \mathcal{I}%
_{0}\left(  M\right)  \right)
\]
is an onto map, where $Aut\left(  \mathcal{I}_{0}\left(  M\right)  \right)  $
is the group of automorphisms of the complex $\mathcal{I}_{0}\left(  M\right)
$ and $\mathcal{MCG}\left(  M\right)  $ is the (extended) mapping class group
of $M,$ i.e. the group of isotopy classes of self-homeomorphisms of $M.$
Moreover, we will show that the map $A_{0}$ has a $\mathbb{Z}_{2}$ kernel.

As explained above, the map
\[
A:\mathcal{MCG}\left(  M\right)  \rightarrow Aut\left(  \mathcal{I}\left(
M\right)  \right)
\]
is not onto. The complex $\mathcal{I}\left(  M\right)  $ is rich enough to
distinguish between an involution and the identity. In other words, the map
$A$ is 1-1. Recall that a finite presentation for $\mathcal{MCG}\left(
M\right)  $ is constructed in \cite{Waj} using the complex of meridians. This
result verifies that the complex of meridians is not sufficient for
characterizing $\mathcal{MCG}\left(  M\right)  .$

For the proof of the above results we list all topological types of surfaces
in the handlebody of genus $2$ and perform a close examination of their links
in $\mathcal{I}\left(  M\right)  .$ This examination establishes that an
automorphism $f$ of $\mathcal{I}\left(  M\right)  $ must map each vertex $v$
in $\mathcal{I}\left(  M\right)  $ to a vertex $f\left(  v\right)  $
consisting of surfaces of the same topological type as those in $v.$ In
particular, $f$ induces an automorphism of the subcomplex $\mathcal{C}\left(
\partial M\right)  $ which permits the use of the corresponding result for
surfaces (see \cite{I1}, \cite{I2} ,\cite{I3}, \cite{Luo}). In the case of 
genus $\geq2$ this analysis cannot
be carried out. However, it can be useful as a basis for establishing
induction in order to show that the automorphisms of the complex
$\mathcal{I}\left(  M\right)  $ are all geometric.

It is well known that for genus $\geq2$ the complex of curves $\mathcal{C}%
\left(  \partial M\right)  $ is a $\delta-$hyperbolic metric space in the
sense of Gromov (see \cite{M-M2},\cite{B}). From the construction of
$\mathcal{I}\left(  M\right)  $ we deduce easily that the embedding of
$\mathcal{C}\left(  \partial M\right)  $ in $\mathcal{I}\left(  M\right)  $ is
isometric and $\mathcal{I}\left(  M\right)  $ is within bounded distance from
the image of $\mathcal{C}\left(  \partial M\right)  $ in $\mathcal{I}\left(
M\right)  ,$ namely, $\mathcal{D}\left(  M\right)  \cup \mathcal{A}\left(
M\right)  .$ In consequence, the complex $\mathcal{I}\left(  M\right)  $ is
itself a $\delta-$hyperbolic metric space in the sense of Gromov. Moreover, it
can be seen (in similar manner as in \cite[Proposition 4.6]{M-M2}) that
$Aut\left(  \mathcal{I}\left(  M\right)  \right)  $ does not contain parabolic
elements but we do not discuss these issues here.

\subsubsection{Notation and terminology}

Let $M$ denote a $3-$dimensional handlebody of genus $n.$ $M$ can be
represented as the union of a handle of index $0$ (i.e. a $3-$ball) with $n$
handles of index $1$ (i.e. $n$ copies of $D^{2}\times \left[  0,1\right]  $).
We fix these handles along with (following the standard terminology) a
longitude $f:S^{1}\rightarrow \partial M$ and meridian $m:S^{1}\rightarrow
\partial M$ for each handle.

For an essential simple closed curve $\alpha$ in $\partial M$ we will be
writing $\left[  \alpha \right]  $ for its isotopy class and the corresponding
vertex in $\mathcal{C}\left(  \partial M\right)  .$ We will be writing
$\left[  S_{\alpha}\right]  $ for the corresponding vertex in $\mathcal{A}%
\left(  M\right)  $ where $S_{\alpha}$ is the annulus corresponding to the
curve $\alpha,$ provided that $\alpha$ is not a meridian boundary. If $\alpha$
is a meridian boundary we will be writing $\left[  D_{\alpha}\right]  $ for
the corresponding vertex in $\mathcal{D}\left(  M\right)  .$ Such curves are
also called meridian curves. By writing $\left[  \alpha \right]  \cap \left[
\beta \right]  =\varnothing$ for non-isotopic curves $\alpha,\beta$ we mean
that there exists curves $\alpha^{\prime},\beta^{\prime}$ isotopic to
$\alpha,\beta$ respectively so that $\alpha^{\prime}\cap \beta^{\prime
}=\varnothing.$ By saying that the class $\left[  \alpha \right]  $ intersects
the class $\left[  \beta \right]  $ at one point we mean that, in addition to
$\left[  \alpha \right]  \cap \left[  \beta \right]  \neq \varnothing,$ there
exist curves $\alpha^{\prime}\in \left[  a\right]  $ and $\beta^{\prime}%
\in \left[  \beta \right]  $ which intersect at exactly one point.

The above notation with square brackets will be similarly used for surfaces.
If $S$ is an incompressible surface we will denote by $Lk\left(  \left[
S\right]  \right)  $ the link of the vertex $\left[  S\right]  $ in
$\mathcal{I}\left(  M\right)  ,$ namely, for each simplex $\sigma$ containing
$\left[  S\right]  $ consider the faces of $\sigma$ not containing $\left[
S\right]  $ and take the union over all such $\sigma.$ We will use the
notation $\ncong$ to declare that two links are not isomorphic as complexes.

As mentioned above, for the rest of this paper, a properly embedded surface
$S$ in $M$ will always mean that, in addition to the above mentioned
requirements, $S$ is isotopic to a surface $\overline{S}$ embedded in
$\partial M$ unless $S$ is a meridian. This assumption further asserts that
such a surface $S$ satisfies the following property

\begin{description}
\item[(SP)] $S$ separates $M$ into two components and the closure of one of
them, denoted by $\Pi_{S},$ is homeomorphic to a product $S\times \left[
0,1\right]  $ with $S\times \left \{  1\right \}  \equiv \overline{S}$ and
$S\times \left \{  0\right \}  \equiv S.$
\end{description}

Consider a surface $S_{\varepsilon}$ properly embedded in $M$ and arbitrarily
close to $\overline{S}$ so that $S_{\varepsilon}$ and $\overline{S}$ bound a
subset of $M,$ say $\Pi_{S_{\varepsilon}},$ homeomorphic to $\overline
{S}\times \left[  0,1\right]  $ with $\overline{S}\times \left \{  1\right \}
\equiv \overline{S}$ and $\overline{S}\times \left \{  0\right \}  \equiv
S_{\varepsilon}.$ It clear that $S_{\varepsilon}$ is properly isotopic to $S$
and we denote this isotopy by $G_{t},t\in \left[  0,1\right]  .$ By standard
isotopy extension properties, see for example \cite[Theorem 1.3 Ch. 8]{Hir},
this proper isotopy is, in fact, an ambient isotopy%
\[
G:M\times \left[  0,1\right]  \rightarrow M
\]
with $G|_{M\times \left \{  0\right \}  }=id_{M}$ and $G\left(  S_{\varepsilon
}\times \left \{  1\right \}  \right)  =S.$ Hence, we have that $G_{1}\left(
\Pi_{S_{\varepsilon}}\right)  $ is homeomorphic to $S\times \left[  0,1\right]
$ with $S\times \left \{  0\right \}  \equiv \overline{S}$ and $S\times \left \{
1\right \}  =G\left(  S_{\varepsilon}\times \left \{  1\right \}  \right)  \equiv
S.$

\section{Invariance of Subcomplexes for genus $n=2$ \label{invariance}}

\bigskip Unless stated otherwise, we restrict our attention to the case where
$M$ is a genus $2$ handlebody. In this section we will show that every
automorphism of $\mathcal{I}\left(  M\right)  $ must preserve the subcomplexes
$\mathcal{A}\left(  M\right)  $ and $\mathcal{D}\left(  M\right)  .$ Moreover,
we will show that for $\left[  I\right]  \in \mathcal{I}\left(  M\right)  ,$
the topological type of the surface $I$ determines the link of $\left[
I\right]  $ in $\mathcal{I}\left(  M\right)  $ and vice-versa. To do this we
will find topological properties for the link of each topological type of
surfaces capable to distinguish their links.

\begin{proposition}
\label{all_meridians}Let $M$ be the handlebody of genus $2$ and $D$ any
meridian. If $S$ is either, an incompressible surface in $M$ with genus $>0$
or, an annular surface with its boundary components being separating curves in
$\partial M$, then $\left[  \partial D\right]  \cap \left[  \partial S\right]
\neq \varnothing.$\newline In particular, a separating curve $\alpha \in \partial
M$ either, bounds a disk or, intersects all meridians.
\end{proposition}

\begin{proof}
First assume that $S$ has genus $>0$ and $D$ is a meridian with $\partial
D\cap \partial S=\varnothing.$ Then, we may assume that $S,D$ intersect
transversely and, hence, $D\cap S$ consists of circles in the interior of $M.$
By irreducibility, we may alter $S$ so that $D\cap S=\varnothing.$ Cutting $M$
along $D$ we obtain that the surface $S$ is properly embedded and
incompressible in a genus $1$ body. This means that $S\ $is an annulus, a
contradiction.\newline Assume now that $S$ is an annular surface with its
boundary components being separating curves in $\partial M$ and $D$ a meridian
with $\partial D\cap \partial S=\varnothing.$ The latter assumption implies
that $D$ is not separating. Cutting $M$ along $D$ we obtain that $S$ is a
surface properly embedded in a solid torus with each component of $\partial S$
being separating, a contradiction.
\end{proof}

\begin{proposition}
\label{algebraic}Let $\alpha$ be a separating curve in $\partial M$ which does
not bound a disk in $M.$ Then each component of $\partial M\setminus \alpha$
contains infinitely many isotopy classes of simple closed curves $\alpha_{i},$
$i,=1,2,\ldots$ such that for each $i,$ $\left[  \partial D\right]
\cap \left[  \alpha_{i}\right]  \neq \varnothing$ for any meridian $D.$
\end{proposition}

\begin{proof}
If for any simple closed curve $\beta$ with $\alpha \cap \beta=$ $\varnothing$
we have that $\left[  \partial D\right]  \cap \left[  \beta \right]
\neq \varnothing$ for any meridian $D,$ we have nothing to show. Thus, let
$\beta$ be a simple closed curve with $\alpha \cap \beta=$ $\varnothing$ and
$\left[  \partial D\right]  \cap \left[  \beta \right]  =\varnothing$ for some
meridian $D.$ Choose a simple closed curve $\gamma$ on $\partial M$ such that
$\alpha \cap \gamma=$ $\varnothing$ and $\beta,\gamma$ intersect at exactly one
point, say $x_{0}.$ The commutator $\beta \gamma \beta^{-1}\gamma^{-1}$ is
freely isotopic to $\alpha.$ Set $x_{0}$ to be the base point of $\pi
_{1}\left(  M,x_{0}\right)  $ and $\pi_{1}\left(  \partial M,x_{0}\right)  .$
Choose generators for $\pi_{1}\left(  M,x_{0}\right)  $ and complete them to a
generating set for $\pi_{1}\left(  \partial M,x_{0}\right)  .$ In the course
of the proof we will consider the free homotopy class of closed curves of the
form $\beta^{i}\gamma^{j},$ $i,j\geq0.$ There is no ambiguity to consider a
curve of the form $\beta^{i}\gamma^{j}$ as an element of $\pi_{1}\left(
M,x_{0}\right)  .$ We will be writing $\left \vert \beta^{i}\gamma
^{j}\right \vert $ to denote the corresponding element of $\pi_{1}\left(
M,x_{0}\right)  .$

Next, cutting $M$ along $D$ we have that $\left \vert \beta \right \vert
=\xi^{n_{0}}$ for some generator $\xi$ of $\pi_{1}\left(  M,x_{0}\right)  ,$
$n_{0}>0$. This follows from the fact that $\alpha \cap \beta=$ $\varnothing$
and, hence, by Proposition \ref{all_meridians}, $\beta$ does not bound a disk.
Hence, it suffices to find infinitely many simple closed curves such that the
corresponding element in $\pi_{1}\left(  M,x_{0}\right)  $ is not the power of
any generator. It is clear that $\left \vert \gamma \right \vert \neq \xi^{m}$ for
all $m\in \mathbb{N}$, otherwise $\left \vert \alpha \right \vert $ would be
trivial in $\pi_{1}\left(  M,x_{0}\right)  $. Choose $\eta$ such that
$\xi,\eta$ generate $\pi_{1}\left(  M,x_{0}\right)  $. Then $\left \vert
\gamma \right \vert $ is a word $w\left(  \xi,\eta \right)  $ containing both
$\xi$ and $\eta.$ Moreover, $\eta$ (and/or $\eta^{-1}$) appears at least twice
in the (reduced) word $\left \vert \gamma^{2}\right \vert .$ To see this, if
$\left \vert \gamma \right \vert $ is a word $w\left(  \xi,\eta \right)  $ in the
letters $\xi,\eta,$ we may write $w\left(  \xi,\eta \right)  $ in the form
\[
w\left(  \xi,\eta \right)  =u\left(  \xi,\eta \right)  \mathrm{\ }v\left(
\xi,\eta \right)  \mathrm{\ }u^{-1}\left(  \xi,\eta \right)
\]
where $v\left(  \xi,\eta \right)  $ is a cyclically reduced word ($u\left(
\xi,\eta \right)  $ being possibly empty). Then
\[
\left \vert \gamma^{2}\right \vert =w\left(  \xi,\eta \right)  \mathrm{\ }%
w\left(  \xi,\eta \right)  =u\left(  \xi,\eta \right)  \mathrm{\ }v^{2}\left(
\xi,\eta \right)  \mathrm{\ }u^{-1}\left(  \xi,\eta \right)  .
\]
Since $\left \vert \gamma \right \vert =w\left(  \xi,\eta \right)  $ contains the
letter $\eta$ (or $\eta^{-1}$) it follows that either $u\left(  \xi
,\eta \right)  $ or, $v\left(  \xi,\eta \right)  $ (or both) must contain $\eta$
(or $\eta^{-1}$). In both cases, $\left \vert \gamma^{2}\right \vert $ contains
$\eta$ and/or $\eta^{-1}$ twice.

To complete the proof we will show that $\left \vert \beta^{i}\gamma
^{2}\right \vert $ is not the power of any generator of $\pi_{1}\left(
M,x_{0}\right)  $ for infinitely many $i$'s. In fact we will restrict
ourselves to odd $i$'s because we need the curve $\beta^{i}\gamma^{2}$ to be simple.

Suppose that $\left \vert \beta^{i}\gamma^{2}\right \vert =\xi_{1}^{n_{1}}$ for
some generator $\xi_{1}=w_{1}\left(  \xi,\eta \right)  $ and $n_{1}\geq1.$
Write $w_{1}\left(  \xi,\eta \right)  $ in the form
\[
u_{1}\left(  \xi,\eta \right)  \mathrm{\ }v_{1}\left(  \xi,\eta \right)
\mathrm{\ }u_{1}^{-1}\left(  \xi,\eta \right)
\]
where $v_{1}\left(  \xi,\eta \right)  $ is cyclically reduced word ($u_{1}$
being possibly empty).

We first examine the case $n_{1}\geq2.$ Let $\ell_{\gamma}$ be the length of
$\left \vert \gamma \right \vert =w\left(  \xi,\eta \right)  .$ For all $i$ such
that $i\cdot n_{0}>6\ell_{\gamma},$ the (reduced) word $\left \vert \beta
^{i}\gamma^{2}\right \vert =\xi^{i\cdot n_{0}}w\left(  \xi,\eta \right)
w\left(  \xi,\eta \right)  $ will have the form
\[
\xi^{N}w^{\prime}\left(  \xi,\eta \right)
\]
where $N>$ $4\ell_{\gamma}$ and $w^{\prime}\left(  \xi,\eta \right)  $ is a
(reduced) word of length $\leq2\ell_{\gamma}$ containing the letter $\eta$
(and/or $\eta^{-1}$) at least twice. Then, the assumption $\left \vert
\beta^{i}\gamma^{2}\right \vert =\xi_{1}^{n_{1}}$ implies that
\[
\xi^{N}w^{\prime}\left(  \xi,\eta \right)  =u_{1}\left(  \xi,\eta \right)
\mathrm{\ }v_{1}^{n_{1}}\left(  \xi,\eta \right)  \mathrm{\ }u_{1}^{-1}\left(
\xi,\eta \right)
\]
with both words being reduced. This implies that $v_{1}\left(  \xi
,\eta \right)  \mathrm{\ }u_{1}\left(  \xi,\eta \right)  $ is a power of $\xi$
and, thus, so is $w_{1}\left(  \xi,\eta \right)  = \xi_{1} .$ It follows that
$\xi_{1}^{n_{1}}=$ $\left \vert \beta^{i}\gamma^{2}\right \vert $ is a power of
$\xi$ and, hence, so is $\left \vert \gamma^{2}\right \vert $, a contradiction.

To complete the proof we need to show that for all (odd) $i$'s considered
above, $\left \vert \beta^{i}\gamma^{2}\right \vert \neq \xi_{1}$ for any
generator $\xi_{1}.$ Assuming that the word $\xi^{N}w^{\prime}\left(  \xi
,\eta \right)  $ is a generator, say $\xi_{1},$ there must exist an
automorphism of the free group $\left \langle \xi,\eta \right \rangle $ mapping
$\xi^{N}w^{\prime}\left(  \xi,\eta \right)  $ onto $\xi.$ Such an automorphism
can be expressed as a product of permutations and Whitehead automorphisms of
$\left \langle \xi,\eta \right \rangle $ (see \cite[pp. 48]{CGKZ}). In other
words, there exists a sequence $\psi_{i},i=1,2,\ldots k$ of permutations and
Whitehead automorphisms such that
\[
\psi_{k}\psi_{k-1}\ldots \psi_{1}\left(  \xi^{N}w^{\prime}\left(  \xi
,\eta \right)  \right)  =\xi.
\]
Recall that a non trivial Whitehead automorphism of the free group
$\left \langle \xi,\eta \right \rangle $ must fix one of the generators, say
$\xi,$ and map $\eta$ to $\eta \xi^{\pm1},$ $\xi^{\pm1}\eta$ or $\xi^{\mp1}%
\eta \xi^{\pm1}.$ To complete the proof of the Proposition, we need the
following three properties where $W\left(  \xi,\eta \right)  $ denotes a
(reduced) word of the form $W\left(  \xi,\eta \right)  =\xi^{M}W^{\prime
}\left(  \xi,\eta \right)  $ where $M>$ $2\ell_{\gamma}$ and $W^{\prime}\left(
\xi,\eta \right)  $ is a (reduced) word of length $\leq2\ell_{\gamma}$
containing the letter $\eta$ (and/or $\eta^{-1}$) at least twice.

\begin{itemize}
\item[(W1)] If $\psi$ is a permutation then the length of $\psi \left(
W\left(  \xi,\eta \right)  \right)  $ is equal to the length of $W\left(
\xi,\eta \right)  .$

\item[(W2)] If $\psi$ is a Whitehead automorphism of $\left \langle \xi
,\eta \right \rangle $ fixing $\eta,$ then the length of $\psi \left(  W\left(
\xi,\eta \right)  \right)  $ is $\geq$ than the length of $W\left(  \xi
,\eta \right)  .$

\item[(W3)] If $\psi$ is a Whitehead automorphism of $\left \langle \xi
,\eta \right \rangle $ fixing $\xi,$ then if $W\left(  \xi,\eta \right)  $
contains a subword of the form $\eta \xi^{m}\eta$ for some $m\in \mathbb{Z},$ so
does the image $\psi \left(  W\left(  \xi,\eta \right)  \right)  .$ \newline If
$W\left(  \xi,\eta \right)  $ contains a subword of the form $\eta \xi^{m}%
\eta^{-1}$ for some $m\in \mathbb{Z}\setminus \left \{  0\right \}  $ then
$\psi \left(  W\left(  \xi,\eta \right)  \right)  $ also contains $\eta \xi
^{m}\eta^{-1}.$
\end{itemize}

It is immediate to check Properties (W1) and (W2). Property (W3) is checked
case by case and it is straightforward.

Since the sequence of automorphisms $\psi_{i},i=1,2,\ldots k$ can be chosen so
that the length of $\xi^{N}w^{\prime}\left(  \xi,\eta \right)  $ decreases with
each application of $\psi_{1},$ $\psi_{2},$ ... we may assume, by property
(W1), that $\psi_{1}$ is not a permutation. Since, for all $i$ such that
$i\cdot n_{0}>6\ell_{\gamma}$ we have that $N>4\ell_{\gamma}$ and $w^{\prime
}\left(  \xi,\eta \right)  $ is a (reduced) word of length $\leq2\ell_{\gamma
},$ it follows by (W2) that $\psi_{1}$ does not fix $\eta.$ Thus $\psi_{1}$
fixes $\xi.$ Let $\lambda \geq1$ be the positive integer such that each
$\psi_{i},i=1,\ldots \lambda$ fixes $\xi$ and $\psi_{\lambda+1}$ does not fix
$\xi.$

By applying the automorphism $\psi_{1}$ on the word $\xi^{N}w^{\prime}\left(
\xi,\eta \right)  $ so that its length strictly reduces we obtain a word of the
form
\[
\xi^{N-1}w^{\prime \prime}\left(  \xi,\eta \right)
\]
where $N$ is assumed to be positive (we work similarly if $N$ is negative).
Recall that $w^{\prime}\left(  \xi,\eta \right)  $ contains the letter $\eta$
(and/or $\eta^{-1}$) at least twice and by (W3), so does $w^{\prime \prime
}\left(  \xi,\eta \right)  .$ Then, the maximum number of consecutive
applications of automorphisms fixing $\xi$ is bounded by the largest power of
$\xi$ in $w^{\prime}\left(  \xi,\eta \right)  .$ In particular, $\lambda$ is
bounded by the length of $w^{\prime}\left(  \xi,\eta \right)  $ which is
$\leq2\ell_{\gamma}$. It follows that the image $\psi_{\lambda}\ldots \psi
_{1}\left(  \xi^{N}w^{\prime}\left(  \xi,\eta \right)  \right)  $ has the
(reduced) form
\[
\xi^{N_{\lambda}}w_{\lambda}^{\prime}\left(  \xi,\eta \right)
\]
where $N_{\lambda}>4\ell_{\gamma}-2\ell_{\gamma}=2\ell_{\gamma}$ and
$w_{\lambda}^{\prime}\left(  \xi,\eta \right)  $ is a (reduced) word of length
$\leq2\ell_{\gamma}$ containing the letter $\eta$ (and/or $\eta^{-1}$) at
least twice. Clearly, $\psi_{\lambda+1}$ must reduce the length of
$\psi_{\lambda}\ldots \psi_{1}\left(  \xi^{N}w^{\prime}\left(  \xi,\eta \right)
\right)  ,$ hence, $\psi_{\lambda+1}$ is not a permutation. By assumption,
$\psi_{\lambda+1}$ does not fix $\xi$, hence, $\psi_{\lambda+1}$ must fix
$\eta.$ Property (W2) gives a contradiction$.$
\end{proof}

\begin{proposition}
\label{non_iso}Let $\left[  D\right]  ,\left[  D^{\prime}\right]
\in \mathcal{D}\left(  M\right)  $, $\left[  S_{\alpha}\right]  ,\left[
S_{\alpha^{\prime}}\right]  \in \mathcal{A}\left(  M\right)  $ and $\left[
I\right]  \in \mathcal{I}\left(  M\right)  \setminus \mathcal{C}\left(  \partial
M\right)  $ where $D$ is a separating meridian, $D^{\prime}$ a non-separating
meridian, $\alpha$ a (non-meridian) separating curve and $\alpha^{\prime}$ a
(non-meridian) non-separating curve. Then the links $Lk\left(  \left[
D\right]  \right)  ,Lk\left(  \left[  D^{\prime}\right]  \right)  ,Lk\left(
\left[  S_{\alpha}\right]  \right)  $, $Lk\left(  \left[  S_{\alpha^{\prime}%
}\right]  \right)  $ and $Lk\left(  \left[  I\right]  \right)  $ are pair-wise
non-isomorphic as complexes.
\end{proposition}

The proof of this proposition is postponed until the end of this Section. As
we allow only surfaces $S$ which can be isotoped to the boundary $\partial M,$
we have four types of surfaces, based on topological type, embedded in
$\partial M$ and, in consequence, four types of vertices in $\mathcal{I}%
\left(  M\right)  \setminus \mathcal{C}\left(  \partial M\right)  :$

\begin{description}
\item[$\left(  T\right)  $] A genus $1$ torus with one boundary component
which is a separating curve in $\partial M.$

\item[$\left(  \Sigma \right)  $] A genus $1$ torus with two mutually isotopic,
non-separating boundary components.

\item[$\left(  P\right)  $] A pair of pants with one boundary component being
a separating curve in $\partial M$ and the other two non-separating and
mutually isotopic.\newline A pair of pants with all $3$ boundary components
being non-separating curves in $\partial M$ and mutually non-isotopic. Such a
pair of pants will be denoted by $\left(  P_{3}\right)  .$

\item[$\left(  Q\right)  $] A sphere with $4$ holes with two pairs of mutually
isotopic, non-separating boundary components.
\end{description}

By analyzing the link of a vertex in $\mathcal{I}\left(  M\right)
\setminus \mathcal{C}\left(  \partial M\right)  $ for each one of the above
types of surfaces we will find that if $I,I^{^{\prime}}$ are surfaces of
different type then $Lk\left(  \left[  I\right]  \right)  \ncong Lk\left(
\left[  I^{^{\prime}}\right]  \right)  .$ In other words, the topological type
of the surface $I,$ when $\left[  I\right]  \in \mathcal{I}\left(  M\right)
\setminus \mathcal{C}\left(  \partial M\right)  ,$ determines the link of
$\left[  I\right]  $ in $\mathcal{I}\left(  M\right)  \setminus \mathcal{C}%
\left(  \partial M\right)  $ and vice-versa (see Corollary \ref{non_iso_more}).

\begin{remark}
\label{r4}Let $\alpha$ be a simple closed curve separating $\partial M$ into
components $T_{\alpha,+},T_{\alpha,-}.$ \newline(a) If $T_{\alpha,+}%
,T\alpha,_{-}$ are isotopic then $M$ is homeomorphic to $T_{\alpha,+}%
\times \left[  0,1\right]  .$\newline This is a well known fact which can be
seen by cutting $M$ along appropriate meridians. \newline(b) The surfaces
$T_{\alpha,+},T_{\alpha,-}$ may or may not be isotopic. \newline To see that
$T_{\alpha,+},T_{\alpha,-}$ may be isotopic, view the handlebody $M$ as the
product $W\times \left[  0,1\right]  ,$ where $W$ is a genus one torus with one
boundary component, and choose $\alpha$ to be a simple closed curve in
$\partial W\times \left[  0,1\right]  .$ Then $\alpha$ bounds a genus $1$
incompressible surface in $M$ and $\alpha$ separates $\partial M$ into
mutually isotopic components $T_{\alpha,+},T_{\alpha,-}.$ \newline To see that
$T_{\alpha,+},T_{\alpha,-}$ may not be isotopic, consider simple closed curves
$\beta,\gamma$ in $T_{\alpha,+}$ such that the commutator $\beta \gamma
\beta^{-1}\gamma^{-1}$ is isotopic to $\alpha$ and $\beta,\gamma$ generate
$\pi_{1}(T_{\alpha,+}).$ If $\beta$ is not a generator for $\pi_{1}(M)$ then
$T_{\alpha,+}$ cannot be isotopic to $T_{\alpha,-}:$ if they were, $M$ would
be homeomorphic to $T_{\alpha,+}\times \left[  0,1\right]  $ and, thus, $\beta$
would have to be a generator of $\pi_{1}(M),$ a contradiction. In the case
$\beta$ is a generator for $\pi_{1}(M)$ we work similarly with the curve
$\beta^{2}$ which is not a generator for $\pi_{1}(M).$
\end{remark}

\begin{remark}
\label{r5}If $S$ is of type $\left(  P_{3}\right)  ,$ then $\partial S$
separates $\partial M\ $into two components which we denote by $P^{+},P^{-}.$
As above, if $P^{+},P^{-}$ are isotopic then $M$ is homeomorphic to the
product $P^{+}\times \left[  0,1\right]  .$ On the other hand, observe that
$P^{+},P^{-}$ may not be isotopic. To see this, choose a non-separating curve
$\alpha$ so that $\alpha$ intersects all meridians (this can be done by
Proposition \ref{algebraic}) and then choose essential non-separating curves
$\beta,\gamma$ such that $\alpha,\beta,\gamma$ are mutually disjoint and non
isotopic. Remove $\alpha,\beta,\gamma$ form $\partial M$ and denote by
$P^{+},P^{-}$ the closures of the two components$.$ Apparently, $\alpha
,\beta,\gamma$ constitute the common boundary of $P^{+},P^{-}.$ If
$P^{+},P^{-}$ were isotopic, then $M$ would be homeomorphic to the product
$P^{+}\times \left[  0,1\right]  $ in which case a meridian not intersecting
$\alpha$ can be found, a contradiction.\smallskip
\end{remark}

We will use the notation $T_{\alpha},$ $\Sigma_{\alpha},$ $P_{\alpha,\beta},$
$P_{\alpha,\beta,\gamma},$ $Q_{\alpha,\beta}$ in order to specify surfaces by
means of their boundary components as follows:

\begin{itemize}
\item Let $\alpha$ be a separating curve in $\partial M.$ View $S_{\alpha}$ as
an annular surface in $\partial M$ and consider the closures of the two
components of $\partial M\setminus S_{\alpha}.$ Each of them is a genus $1$
torus with one boundary component isotopic to $\alpha.$ We will be denoting
them by $T_{\alpha,+}$ and $T_{\alpha,-}.$ Note that $T_{\alpha,+}$ may or,
may not be isotopic to $T_{\alpha,-}.$

\item Let $\alpha$ be a non separating curve in $\partial M.$ Similarly, the
closure of $\partial M\setminus S_{\alpha}$ is a genus $1$ torus with two
boundary components both isotopic to $\alpha.$ We will be denoting this
surface by $\Sigma_{\alpha}.$

\item Let $\alpha$ be a separating curve and $\beta$ a non separating curve in
$\partial M$ with $\alpha \cap \beta=\varnothing.$ If $T_{\alpha,+}$ is the
subsurface of $\partial M$ containing $\beta,$ view $S_{\beta}$ as an annular
surface in $T_{\alpha,+}$ and set $P_{\alpha,\beta}$ to be the closure of
$T_{\alpha,+}\setminus S_{\beta}.$ Apparently, $P_{\alpha,\beta}$ is a pair of
pants with two boundary components isotopic to $\beta$ and the third isotopic
to $\alpha.$

\item Let $\alpha,\beta,\gamma$ be non separating and mutually non isotopic,
disjoint curves in $\partial M.$ View $S_{\alpha},$ $S_{\beta},$ $S_{\gamma}$
as annular surfaces in $\partial M$ and consider the closures of the two
components of $\partial M\setminus \left(  S_{\alpha}\cup S_{\beta}\cup
S_{\gamma}\right)  .$ Each of them is a pair of pants with boundary components
isotopic to $\alpha,\beta,\gamma.$ We will be denoting them by $P_{\alpha
,\beta,\gamma}^{+}$ and $P_{\alpha,\beta,\gamma}^{-}.$ Note that
$P_{\alpha,\beta,\gamma}^{+}$ may or, may not be isotopic to $P_{\alpha
,\beta,\gamma}^{-}.$

\item Let $\alpha,\beta$ be non separating, disjoint and non isotopic curves
in $\partial M.$ View $S_{\alpha},$ $S_{\beta}$ as annular surfaces in
$\partial M$ and set $Q_{\alpha,\beta}$ to be the closure of $\partial
M\setminus \left(  S_{\alpha}\cup S_{\beta}\right)  .$ Apparently,
$Q_{\alpha,\beta}$ is a surface of genus $0$ with four boundary components
\end{itemize}

In the following subsections we will be using repeatedly the following Lemma.
Recall that if $S\ $is not a meridian then there exists $\overline{S}$
isotopic to $S$ with $\overline{S}$ being embedded in $\partial M$ so that
$S\cup \overline{S}$ bounds a subset $\Pi_{S}$ of $M$ homeomorphic to the
product $S\times \left[  0,1\right]  $ (see property (SP) above).

\begin{lemma}
\label{sigma_torus}Let $S,S^{\prime}$ be properly embedded incompressible
surfaces in $M$ such that none of them is a meridian. Then $\left[  S\right]
\cap \left[  S^{\prime}\right]  =\varnothing$ if and only if there exists
$\overline{S}\subset \partial M$ isotopic to $S$ and $\overline{S}^{\prime
}\subset \partial M$ isotopic to $S^{\prime}$ such that either, $\overline
{S}\subset \overline{S}^{\prime}$ or, $\overline{S}^{\prime}\subset \overline
{S}$ or, $\overline{S}\cap \overline{S}^{\prime}=\varnothing.$
\end{lemma}

\begin{proof}
For the "if" part, consider the product sets $\Pi_{S}=S\times \left[
0,1\right]  $ and $\Pi_{S^{\prime}}=S^{\prime}\times \left[  0,1\right]  .$ It
is clear that the condition $\overline{S}\cap \overline{S}^{\prime}%
=\varnothing$ implies that up to isotopy, $\Pi_{S}$ and $\Pi_{S^{\prime}}$ are
disjoint. Thus, $S\times \left \{  0\right \}  \equiv S$ and $S^{\prime}%
\times \left \{  0\right \}  \equiv S^{\prime}$ are disjoint as required.
Similarly, if $\overline{S}\subset \overline{S}^{\prime}$ then, up to isotopy,
$\Pi_{S}\subset S^{\prime}\times \left(  0,1\right]  $ and, hence,
$S\times \left \{  0\right \}  \equiv S$ is disjoint from $S^{\prime}%
\times \left \{  0\right \}  \equiv S^{\prime}$

For the "only if" part first observe that given a properly embedded
incompressible (non-meridian) surface $S$ in $M$ then for any choice of
$\overline{S}\in \left[  S\right]  ,$ the surfaces $\overline{S}$ and the
closure of $\partial M\setminus \overline{S}$ are of distinct topological type,
except in the following two cases

\begin{description}
\item[Case A:] $S$ is of type $\left(  P_{3}\right)  $, i.e. all $3$ boundary
components are non-separating curves in $\partial M$.

\item[Case B:] $S$ is of type $(T).$
\end{description}

For, if $S$ is of type $(P)$ but not of type $\left(  P_{3}\right)  $ then the
closure of $\partial M\setminus \overline{S}$ is a surface with two components
of genus $1$ and $0$; if $S$ is of type $(Q)$ then the closure of $\partial
M\setminus \overline{S}$ has two components both of genus $0.$ In the case $S$
is of type $(\Sigma)$ (resp. annular surface) the closure of $\partial
M\setminus \overline{S}$ is an annular surface (resp. either of type $(\Sigma)$
or disconnected). Hence, if $S$ is, neither of type $\left(  P_{3}\right)  ,$
nor of type $(T)$ then the boundary components of $S$ along with its
topological type determine uniquely the component of $\partial M\setminus
\partial S$ whose closure is isotopic to $S.$ As $\left[  S\right]
\cap \left[  S^{\prime}\right]  =\varnothing,$ it is clear that $\left[
\partial S\right]  \cap \left[  \partial S^{\prime}\right]  =\varnothing$ and,
therefore, a case by case examination reveals that if none of $S,S^{\prime}$
is of type $(T)$ or $\left(  P_{3}\right)  $ then $\overline{S},\overline
{S}^{\prime}$ are either disjoint or one is contained in the other.\smallskip

\noindent \underline{Case A:} $S$ is of type $\left(  P_{3}\right)  .$ Then
$\partial S$ separates $\partial M\ $into two components which we denote by
$P^{+},P^{-}$ and assume that $\overline{S}=P^{+}.$ As mentioned above,
$P^{+},P^{-}$ may or may not be isotopic.

\begin{itemize}
\item If $S^{\prime}$ is annular then, up to isotopy, either $\overline
{S}^{\prime}\subset P^{+}$ or, $\overline{S}^{\prime}\subset P^{-}.$ In either
case we have $\overline{S}^{\prime}\subset \overline{S}$ or $\overline{S}%
\cap \overline{S}^{\prime}=\varnothing$ as required.

\item If $S^{\prime}$ is of type $\left(  P\right)  $ but not of type $\left(
P_{3}\right)  $ then the separating boundary curve of $S^{\prime}$ will
necessarily intersect $\partial S$ violating the assumption $\partial
S\cap \partial S^{\prime}=\varnothing.$

\item If $S^{\prime}$ is of type $\left(  T\right)  $ then the separating
curve $\partial S^{\prime}$ will necessarily intersect $\partial S,$ a contradiction.

\item If $S^{\prime}$ is of type $\left(  Q\right)  $ then assumption
$\partial S\cap \partial S^{\prime}=\varnothing$ implies that one component of
$\partial S$ will be in the interior of $\overline{S}^{\prime}$ and the other
two isotopic to the the boundary components of $S^{\prime}$ respectively.
Hence, all 3 components of $\partial S$ are, up to isotopy, contained in
$\overline{S}^{\prime}.$ Apparently, $\overline{S}\subset \overline{S}^{\prime
}$ as required.\newline

\item If $S^{\prime}$ is of type $\left(  \Sigma \right)  $ we have, in a
similar manner, that $\overline{S}\subset \overline{S}^{\prime}.$
\end{itemize}

\noindent \underline{Case B:} $S$ is of type $(T).$ In this case $\partial S$
separates $\partial M\ $into $2$ boundary components both being of type $(T).$
Denote them by $T_{+},T_{-}$ and, as mentioned above, $T_{+},T_{-}$ may or may
not be isotopic. By property (SP), $S$ is isotopic to at least one of
$T_{+},T_{-}.$ By changing notation, if necessary, we assume that $T_{+}%
\in \left[  S\right]  $ (and then $T_{-}$ may or may not belong to $\left[
S\right]  $) and $S,T_{+}$ bound a product $S\times \left[  0,1\right]
\equiv \Pi_{S}.$ Observe that if

\begin{itemize}
\item $S^{\prime}$ is annular then, apparently, $\overline{S}^{\prime}%
\subset \overline{S}$ or, $\overline{S}\cap \overline{S}^{\prime}=\varnothing.$

\item $S^{\prime}$ is of type $\left(  P\right)  $ (note that since $\partial
S\cap \partial S^{\prime}=\varnothing,$ $S^{\prime}$ cannot be of type $\left.
\left(  P_{3}\right)  \right)  ,$ then the $2$ non-separating (and mutually
isotopic) components of $S^{\prime}$ belong to either $T_{+},$ in which case
$\overline{S}^{\prime}\subset \overline{S}$ or, to $T_{-},$ in which case
$\overline{S}\cap \overline{S}^{\prime}=\varnothing.$

\item $S^{\prime}$ is of type $\left(  Q\right)  $ then $S^{\prime}$
intersects both $T_{+}$ and $T_{-}$ and, therefore, $S^{\prime}$ intersects
$S\times \left \{  0\right \}  \equiv S$ a contradiction.\newline
\end{itemize}

We conclude the proof of Case B (and, hence, the lemma) by examining the case
where $S$ is of type $(T)$ and $S^{\prime}$ is of type $(\Sigma).$ $S^{\prime
}$ has two boundary components, hence, $\overline{S}^{\prime}$ cannot be
contained in $T_{+}\equiv \overline{S}$ which has one boundary component.

If $\overline{S}\cap \overline{S}^{\prime}=\varnothing$ we are done. Assume
that $\overline{S}\cap \overline{S}^{\prime}\neq \varnothing$ and we will show
that $\overline{S}\subset \overline{S}^{\prime}.$ If $\partial \overline
{S}^{\prime}\subset T_{-}$ then $\overline{S}^{\prime}\supset T_{+}%
\equiv \overline{S}$ as required. Assume that $\partial \overline{S}^{\prime
}\subset T_{+}.$ As $\Pi_{S^{\prime}}\equiv S^{\prime}\times \left[
0,1\right]  $ is connected and $\partial S\subset S^{\prime}\times \left \{
1\right \}  ,$ it follows that $S\subset \Pi_{S^{\prime}}.$ Choose separating
simple closed curve $\alpha$ in $\overline{S}^{\prime}\equiv S^{\prime}%
\times \left \{  1\right \}  $ to form a product set denoted by $\alpha
\times \left[  0,1\right]  $ such that $\alpha \times \left \{  1\right \}
\equiv \alpha$ and $\alpha \times \left \{  t\right \}  \subset$ $S^{\prime}%
\times \left \{  t\right \}  $ for all $t\in \left[  0,1\right]  .$ Cut the
product space $\Pi_{S^{\prime}}$ along the annulus $\alpha \times \left[
0,1\right]  $ to obtain a handlebody of genus $2.$ We view this handlebody as
the product $T_{-}\times \left[  0,1\right]  $ and $S$ is properly embedded in
$T_{-}\times \left[  0,1\right]  $ with $\partial S\subset T_{-}\times \left \{
0\right \}  .$ $S\ $is properly isotopic in $T_{-}\times \left[  0,1\right]  $
with a surface $R$ with $\partial R\subset \partial T_{-}\times \left[
0,1\right]  .$

\begin{figure}[ptb]
\begin{center}
\includegraphics[scale=0.6]
{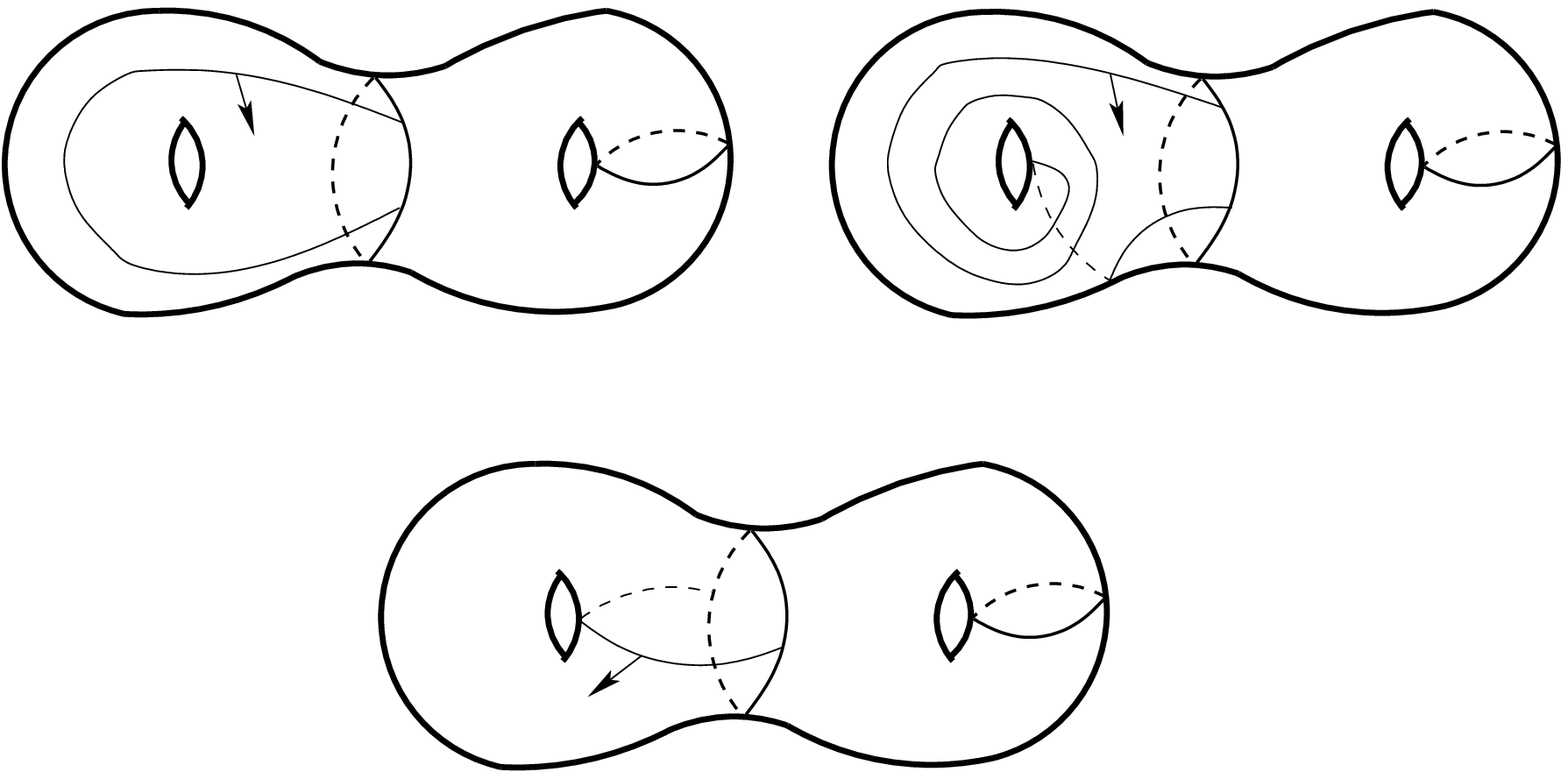}
\end{center}
\par
\begin{picture}(22,12)
\put(100,148){$\sigma$}
\put(120,115){$D_i$}
\put(284,116){$D_i$}
\put(194,27){$D_i$}
\put(197,147){$D$}
\put(362,147){$D$}
\put(272,57){$D$}
\put(270,146){$\sigma$}
\put(159,38){$\tau$}
\put(39,103){(a)}
\put(39,37){(b)}
\end{picture}
\caption{$\,$}%
\label{fig02}%
\end{figure}

Thus, it suffices to show that if $R$ is a genus $1$ surface properly embedded
in a genus $2$ handlebody $T_{-}\times \left[  0,1\right]  $ with $\partial
R\subset \partial T_{-}\times \left[  0,1\right]  $ then $R$ is isotopic to
$T_{-}\times \left \{  0\right \}  .$ Choose standard non-separating, non
properly isotopic meridians $D_{1},D_{2}$ in $T_{-}\times \left[  0,1\right]  $
each intersecting $R$ into two essential arcs, denoted by $\sigma_{1}%
,\sigma_{2}$ with boundary points on $\partial R.$ Moreover, $\left(
D_{1}\cup D_{2}\right)  \cap \partial R$ consists of $4$ points which are
precisely the boundary points of $\sigma_{1},\sigma_{2}.$ These $4$ points
separate $\partial R$ into $4$ subarcs denoted by $\tau_{j},$ $j=1,2,3,4.$ Cut
$T_{-}\times \left[  0,1\right]  $ along $D_{1},D_{2}$ to obtain a $3-$ball
with two copies $D_{1}^{A},D_{1}^{B}$ of $D_{1}$ and two copies $D_{2}%
^{A},D_{2}^{B}$ of $D_{2}$ marked on its boundary. Denote by $R_{\delta}$ the
surface, which is just a disk, in the $3-$ball corresponding to $R.$ Then the
boundary of $R_{\delta}$ is the curve $\delta$ obtained by the juxtaposition
of the arcs%
\[
\delta=\sigma_{1}^{A}\cup \tau_{1}\cup \sigma_{2}^{B}\cup \tau_{3}\cup \sigma
_{1}^{A}\cup \tau_{2}\cup \sigma_{2}^{B}\cup \tau_{4}%
\]
where $\sigma_{1}^{A},\sigma_{1}^{B}$ are the arcs in $D_{1}^{A},D_{1}^{B}$
resp. induced (after cutting) by $\sigma_{1},\sigma_{2}.$ Fix an orientation
transverse to $R.$ This induces an orientation transverse to $R_{\delta}.$
Denote by $\sigma_{1}^{A}(N)$ the subarc of $\partial D_{1}^{A}$ determined by
the orientation of $R_{\delta}$ which, of course, is homotopic to $\sigma
_{1}^{A}$ with endpoints fixed. Similarly for $\sigma_{1}^{B}(N),$ $\sigma
_{2}^{A}(N)$ and $\sigma_{2}^{B}(N).$ As $R$ is orientable, all these $4$
subarcs are contained in one of the two hemispheres of the boundary of the
$3-$ball determined by $\partial R_{\delta}.$ Denote this hemisphere by $N.$
We may isotope $R_{\delta}$ arbitrarily close to the boundary of $N$ to obtain
a disk $R_{\delta \left(  N\right)  }$ whose boundary is the curve
\[
\delta \left(  N\right)  =\sigma_{1}^{A}\left(  N\right)  \cup \tau_{1}%
\cup \sigma_{2}^{B}\left(  N\right)  \cup \tau_{3}\cup \sigma_{1}^{A}\left(
N\right)  \cup \tau_{2}\cup \sigma_{2}^{B}\left(  N\right)  \cup \tau_{4}.
\]
\newline Apparently, $R_{\delta \left(  N\right)  }$ is isotopic (with boundary
fixed) to the disk in $N$ bounded by $\delta \left(  N\right)  .$ After gluing
back $D_{1}^{A},D_{1}^{B}$ and $D_{2}^{A},D_{2}^{B}$ the disk in $N$ bounded
by $\delta \left(  N\right)  $ is the union of $T_{-}$ along with an annulus in
$\partial T_{-}\times \left[  0,1\right]  $ bounded by $\partial R$ and
$\partial T_{-}\times \left \{  0\right \}  .$ Hence, $R\ $is isotopic to $T_{-}$
as required.
\end{proof}

\subsection{Separating Meridians\label{sm}}

Let $D$ be a separating meridian in $M.$ We will study the link of the vertex
$\left[  D\right]  $ in $\mathcal{I}\left(  M\right)  .$ $D$ decomposes $M$
into two solid tori $T^{+},T^{-}.$ Denote by $D^{+}$ (resp. $D^{-}$) the
unique meridian in $T^{+}$ (resp. $T^{-}$). Consider the infinite sequence
$\alpha_{i}^{+},i=0,1,2,\ldots$ (resp. $\alpha_{j}^{-},j=0,1,2,\ldots$)
consisting of all essential simple closed and mutually non-isotopic curves in
$\partial T^{+}\setminus D$ (resp. $\partial T^{-}\setminus D$) excluding
those isotopic to $\partial D$. For $i=0$ (resp. $j=0$) denote by $\left[
S_{\alpha_{0}^{+}}\right]  $ (resp. $\left.  \left[  S_{\alpha_{0}^{-}%
}\right]  \right)  $ the vertex $\left[  D^{+}\right]  $ (resp. $\left.
\left[  D^{-}\right]  \right)  .$ The corresponding annular surfaces
$S_{\alpha_{i}^{+}},$ $i=1,2,\ldots$ and $S_{\alpha_{j}^{-}}$, $j=1,2,\ldots$
give rise to distinct vertices $\left[  S_{\alpha_{i}^{+}}\right]  ,$ $\left[
S_{\alpha_{j}^{-}}\right]  $ which all belong to $Lk\left(  D\right)  .$ If
$S$ is an incompressible surface in $M$ with $\left[  D\right]  \cap \left[
S\right]  =\varnothing$ then $S$ is contained in a solid torus hence, $S$ is
either, an annulus or, a meridian. Hence, the vertex set of $Lk\left(  \left[
D\right]  \right)  $ is the following set
\[
\left \{  \left[  S_{\alpha_{i}^{+}}\right]  \bigm \vert i=0,1,2,\ldots
\right \}  \cup \left \{  \left[  S_{\alpha_{j}^{-}}\right]  \bigm \vert
j=0,1,2,\ldots \right \}
\]
Apparently, for any $i,j$ the surfaces $S_{\alpha_{i}^{+}},$ $S_{\alpha
_{j}^{-}}$ are disjoint hence the $Lk(\left[  D\right]  )$ contains all edges
$\left(  \left[  S_{\alpha_{i}^{+}}\right]  ,\left[  S_{\alpha_{j}^{-}%
}\right]  \right)  $ $\forall i,j=1,2,\ldots.$ Moreover, for all
$i,i^{{\prime}}$ with $i\neq i^{{\prime}}$ the surfaces $S_{\alpha_{i}^{+}%
},S_{\alpha_{i^{{\prime}}}^{+}}$ intersect and similarly for $S_{\alpha
_{j}^{-}},S_{\alpha_{j^{{\prime}}}^{-}}.$ Since it is clear that no
$2-$dimensional simplices exist in $Lk\left(  \left[  D\right]  \right)  ,$ we
have shown the following property

\begin{description}
\item[($sM$-$1$)] If $D$ is a separating meridian ($sM$) then $Lk\left(
\left[  D\right]  \right)  $ is isomorphic to the bi-infinite complete
bipartite graph.
\end{description}

\noindent In particular, if the length of a path is given by the number of its
edges, we have

\begin{description}
\item[($sM$-$2$)] $Lk\left(  \left[  D\right]  \right)  $ does not contain
simple closed cycles of length $3.$
\end{description}

\noindent Properties ($sM$-$1$), ($sM$-$2$) will be used later in Section
\ref{invariance_proof} to prove Proposition \ref{non_iso}. Analogous
properties will be stated at the end of each of the following subsections.

\begin{figure}[ptb]
\begin{center}
\includegraphics[scale=0.6]
{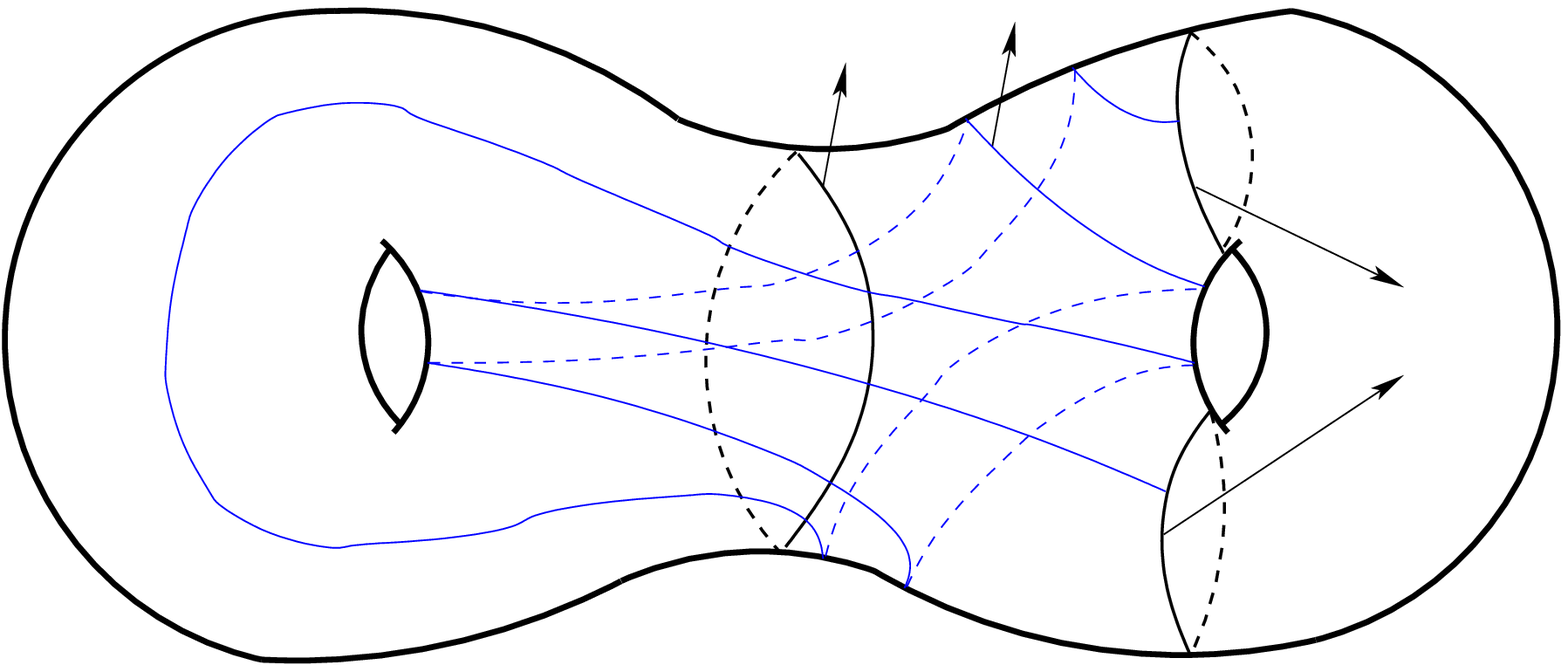}
\end{center}
\par
\begin{picture}(22,12)
\put(213,150){$D_i$}
\put(327,90){$D$}
\put(249,157){$\rho$}
\end{picture}
\caption{$\,$}%
\label{fig03}%
\end{figure}

\subsection{Non-separating Meridians\label{nsm}}

Let $D$ be a non-separating meridian in $M.$ Consider the infinite sequence
$D_{i},$ $i=1,2,\ldots$ consisting of all separating meridians each being
disjoint from $D$ and having pair-wise non-isotopic boundaries. Each $D_{i}$
separates $M$ into two solid tori $T_{i,+},T_{i,-}$ and we may assume that $D$
is a meridian in $T_{i,-}.$ As before, let $\alpha_{j}^{i},j=0,1,2,\ldots$ be
the infinite sequence of essential simple closed curves in $\partial
T_{i,+}\setminus D_{i}$ which are not isotopic to $\partial D_{i}$ and let
$S_{\alpha_{j}^{i}}$, $j=1,2,\ldots$ be the corresponding incompressible
annular surfaces in $T_{i,+}.$ For $j=0$ we set $S_{\alpha_{0}^{i}}$ to be the
unique (up to isotopy) meridian in $T_{i,+}.$

Let now $S$ be an incompressible surface in $M$ which is not isotopic to any
$D_{i},i=1,2,\ldots.$ Since $\left[  D\right]  \cap \left[  S\right]
=\varnothing,$ $S$ is contained in a solid torus and, hence, $S$ is either, an
annulus or, a meridian (in the solid torus). In case $S$ is an annulus in the
solid torus it must also be an annular surface in $M$ because $S$ is assumed
to be isotopic with a surface in $\partial M.$ Hence, the vertex set of
$Lk\left(  \left[  D\right]  \right)  $ is the following set
\[
\cup_{i=1}^{\infty}\left \{  \left[  D_{i}\right]  ,\left[  S_{\alpha_{0}^{i}%
}\right]  ,\left[  S_{\alpha_{1}^{i}}\right]  ,\left[  S_{\alpha_{2}^{i}%
}\right]  ,\ldots \right \}  .
\]
It is clear that for each $i,i=1,2,\ldots$the $Lk\left(  \left[  D\right]
\right)  $ contains all edges $\left(  \left[  D_{i}\right]  ,\left[
S_{\alpha_{j}^{i}}\right]  \right)  $ $\forall j=0,1,2,\ldots.$ It is easy to
see that $\left[  D_{i}\right]  \cap \left[  D_{i^{{\prime}}}\right]
\neq \varnothing$ if $i\neq i^{{\prime}}.$ Hence, $Lk\left(  \left[  D\right]
\right)  $ does not contain any edge of the form $\left(  \left[
D_{i}\right]  ,\left[  D_{i^{{\prime}}}\right]  \right)  $ for all indexes
$i\neq i^{{\prime}}.$\newline \noindent \textbf{Claim}: for any $i,i^{{\prime}}$
with $i\neq i^{{\prime}}$ there exists at most one pair of indices $j\left(
i\right)  $ and $j\left(  i^{{\prime}}\right)  $ so that $\left[
S_{\alpha_{j\left(  i\right)  }^{i}}\right]  =\left[  S_{\alpha_{j\left(
i^{{\prime}}\right)  }^{i^{{\prime}}}}\right]  .$

\noindent \textbf{Proof of Claim}: First observe that, as $D_{i}\cap
D_{i^{{\prime}}}\neq \varnothing,$ $\partial D_{i^{{\prime}}}$ must contain
either, a path $\sigma$ of the forms shown in Figure \ref{fig02}a or, a path
$\tau$ of the form shown in Figure \ref{fig02}b (otherwise, $D_{i}%
,D_{i^{{\prime}}}$ would belong to the same isotopy class). Note that the
defining difference between a path $\tau$ and a path $\sigma$ is that the
juxtaposition of a path $\tau$ of the form shown in Figure \ref{fig02}b with
one of the two subarcs of $\partial D_{i}$ gives rise to a meridian boundary.
In the case $\partial D_{i^{{\prime}}}$ contains only paths of the form shown
in Figure \ref{fig02}b then $S_{\alpha_{0}^{i}}$ is also a meridian for
$T_{i^{\prime},+}.$ In other words $\left[  S_{\alpha_{0}^{i}}\right]
=\left[  S_{\alpha_{0}^{i^{{\prime}}}}\right]  .$

In the case $\partial D_{i^{{\prime}}}$ contains only paths of the form shown
in Figure \ref{fig02}a then there exists a unique essential curve
$\alpha_{j\left(  i\right)  }^{i}$ in $T_{i,+}$ disjoint from $\partial
D_{i^{{\prime}}}.$ Hence this curve is an essential curve in $T_{i^{\prime}%
,+},$ i.e., $\alpha_{j\left(  i\right)  }^{i}$ is isotopic to $\alpha
_{j\left(  i^{{\prime}}\right)  }^{i}$ for some $j\left(  i^{{\prime}}\right)
.$ In other words, $\left[  S_{\alpha_{j\left(  i\right)  }^{i}}\right]
=\left[  S_{\alpha_{j\left(  i^{{\prime}}\right)  }^{i^{{\prime}}}}\right]  $
as claimed.

In the case $\partial D_{i^{{\prime}}}$ contains both types of paths shown in
Figure \ref{fig02}a, \ref{fig02}b then there is no pair of indexes as in the
statement of the Claim (see, for example, Figure \ref{fig03} where
$D_{i^{{\prime}}}$ is obtained by two copies of $D$ joined by $\rho$). This
completes the proof of the Claim.

We will now state three properties for $Lk\left(  \left[  D\right]  \right)  $
when $D$ is a non-separating meridian ($nsM$) in $M$ to be used later in
Section \ref{invariance_proof} to prove Proposition \ref{non_iso}.

\begin{description}
\item[($nsM$-$1$)] $Lk\left(  \left[  D\right]  \right)  $ is not isomorphic
to the bi-infinite complete bipartite graph.
\end{description}

\noindent This follows easily from the fact that $D_{2}$ intersects
$S_{\alpha_{j}^{1}}$ for infinitely many $j.$ Moreover,

\begin{description}
\item[\textbf{(}$nsM$-$2$\textbf{)}] any simple closed cycle in $Lk\left(
\left[  D\right]  \right)  $ has length at least $4.$
\end{description}

\noindent This follows easily from the fact that $3$ is the maximum number of
mutually disjoint essential simple closed curves to be found on $\partial M$
(the vertices of a closed $3-$cycle along with $\left[  D\right]  $ give a
contradiction). In fact, it can be shown that any simple closed cycle in
$Lk\left(  \left[  D\right]  \right)  $ has length at least $5,$ but we do not
need this.

\begin{description}
\item[($nsM$-$3$)] $Lk\left(  \left[  D\right]  \right)  $ contains infinitely
many vertices of infinite valence.
\end{description}

\subsection{Annular surfaces with separating boundary}

Let $\alpha$ be a separating curve in $\partial M$ which is not homotopically
trivial in $M.$ Let $S_{\alpha}$ be the properly embedded annular
(incompressible) surface whose boundary consists of two (parallel) copies of
$\alpha$ so that $S_{\alpha}$ is isotopic to the annular region in $\partial
M$ bounded by its boundary components. Denote by $T_{\alpha,+},T_{\alpha,-}$
the closures of the components of $\partial M\setminus S_{\alpha}$ (as
explained in the text following Remark \ref{r5}). Denote by $\alpha_{i}%
^{+},i=1,2,\ldots$ (resp. $\alpha_{j}^{-},j=1,2,\ldots$) the infinite sequence
of pair-wise non-isotopic essential simple closed curves in $T_{\alpha,+}$
(resp. $T_{\alpha,-}$) which intersect pair-wise. We will first list all
vertices of the complex $\mathcal{I}\left(  M\right)  $ which belong to
$Lk\left(  \left[  S_{\alpha}\right]  \right)  .$

\begin{itemize}
\item Both components $T_{\alpha,+},T_{\alpha,-}$ can be viewed as properly
embedded surfaces in $M.$ By Proposition \ref{all_meridians}, these surfaces
are incompressible. Hence, they determine vertices $\left[  T_{\alpha
,+}\right]  ,\left[  T_{\alpha,-}\right]  $ in $Lk\left(  \left[  S_{\alpha
}\right]  \right)  .$ Recall that these two vertices are not necessarily
distinct (see Remark \ref{r4}).

\item The annular surfaces $S_{\alpha_{i}^{+}},$ $i=1,2,\ldots$ and
$S_{\alpha_{j}^{-}}$, $j=1,2,\ldots$ give rise to distinct (with respect to
$i,j$) vertices $\left[  S_{\alpha_{i}^{+}}\right]  ,\left[  S_{\alpha_{j}%
^{-}}\right]  .$

\item For each $i=1,2,\ldots$ (resp. $j=1,2,\ldots$) the surfaces
$P_{\alpha,\alpha_{i}^{+}}$ (resp. $P_{\alpha,\alpha_{j}^{-}}$) are, by
Proposition \ref{all_meridians}, incompressible and rise to distinct vertices
$\left[  P_{\alpha,\alpha_{i}^{+}}\right]  $ (resp. $\left.  \left[
P_{\alpha,\alpha_{j}^{-}}\right]  \right)  .$

\item For each $i=1,2,\ldots$ (resp. $j=1,2,\ldots$) the surfaces
$\Sigma_{\alpha_{i}^{+}}$ (resp. $\Sigma_{\alpha_{j}^{-}}$) may or may not be
incompressible depending on whether or not $\alpha_{i}^{+}$ (resp. $\alpha
_{j}^{-})$ intersects all meridians. However, by Proposition \ref{algebraic},
$\Sigma_{\alpha_{i}^{+}}$ (resp. $\Sigma_{\alpha_{j}^{-}}$) is incompressible
for infinitely many $i$ 's (resp. $j$ 's) and, hence, we obtain distinct
vertices $\Sigma_{\alpha_{i}^{+}}$ (resp. $\Sigma_{\alpha_{j}^{-}}$) for
infinitely many $i$'s (resp. $j$ 's).

\item For each $i,j=1,2,\ldots$ the surfaces $Q_{\alpha_{i}^{+},\alpha_{j}%
^{-}}$ may or may not be incompressible depending on whether or not the union
$\alpha_{i}^{+}\cup \alpha_{j}^{-}$ intersects all meridians. By Proposition
\ref{algebraic}, $Q_{\alpha,i,j}$ is incompressible for infinitely many $i$'s
and $j$'s and, hence, we obtain distinct vertices $\left[  Q_{\alpha_{i}%
^{+},\alpha_{j}^{-}}\right]  $ for infinitely many $i$'s and $j$'s.
\end{itemize}

This is a complete list of the vertices in $Lk\left(  S_{\alpha}\right)  .$ To
see this, let $S\ $be a properly embedded incompressible surface so that
$\left[  S\right]  \cap \left[  S_{\alpha}\right]  =\varnothing.$ As
$S_{\alpha}$ is an annular incompressible surface with $\alpha$ separating, by
Proposition \ref{all_meridians}, $S$ is not a meridian. Hence, by definition
of the vertices of the complex $\mathcal{I}\left(  M\right)  ,$ $S$ is
isotopic to a surface embedded in $\partial M.$

If $\partial S$ is connected then $S$ has genus $\geq1$ and, as $\partial S$
must separate $\partial M$, $\partial S$ is isotopic to $\alpha.$ Thus, $S$ is
isotopic to either $T_{\alpha,+}$ or, $T_{\alpha,-}.$

If $\partial S$ has two components then, up to isotopy, they are both
contained in either $T_{\alpha,+}$ or, $T_{\alpha,-}$ and, hence, the two
boundary components are isotopic to either a curve $\alpha_{i}^{+}$ in
$T_{\alpha,+}$ or, a curve $\alpha_{j}^{-}$ in $T_{\alpha,-}.$ If the genus of
$S\ $is $0,$ then $S$ is an annulus in either $T_{\alpha,+}$ or, $T_{\alpha
,-}$ Therefore, $S$ is isotopic to either $S_{\alpha_{i}^{+}}$ or,
$S_{\alpha_{j}^{-}}$ for some $i$ or $j.$ If the genus of $S\ $is $1,$ then
$S$ is isotopic to either $\Sigma_{\alpha_{i}^{+}}$ or, $\Sigma_{\alpha
_{j}^{-}}$ for some $i$ or $j.$

Similarly, if $\partial S$ has 3 components then $S$ is isotopic to either
$P_{\alpha,\alpha_{i}^{+}}$ or $P_{\alpha,\alpha_{j}^{-}}$ for some $i$ or
$j.$

Finally, if $\partial S$ has 4 components then $S$ is isotopic to
$Q_{\alpha_{i}^{+},\alpha_{j}^{-}}$ for some $i,j.$ Note that all
incompressible surfaces considered here cannot have more than $4$ boundary components.

We will not list all edges in $Lk\left(  \left[  S_{\alpha}\right]  \right)
.$ However, it is clear that $Lk\left(  \left[  S_{\alpha}\right]  \right)  $
contains all edges of the form $\left(  \left[  S_{\alpha_{i}^{+}}\right]
,\left[  S_{\alpha_{j}^{-}}\right]  \right)  $ for all $i,j=1,2,\ldots.$
Moreover, for all $i,i^{{\prime}}$ with $i\neq i^{{\prime}}$ we have $\left[
S_{\alpha_{i}^{+}}\right]  \cap \left[  S_{\alpha_{i^{{\prime}}}^{+}}\right]
\neq$ $\varnothing$ and similarly for $\left[  S_{\alpha_{j}^{-}}\right]
,\left[  S_{\alpha_{j^{{\prime}}}^{-}}\right]  .$ Hence, $Lk\left(  \left[
S_{\alpha}\right]  \right)  $ does not contain any edge of the form $\left(
\left[  S_{\alpha_{i}^{+}}\right]  ,\left[  S_{\alpha_{i^{{\prime}}}^{+}%
}\right]  \right)  ,$ $i\neq i^{{\prime}}$ or $\left(  \left[  S_{\alpha
_{j}^{-}}\right]  ,\left[  S_{\alpha_{j^{{\prime}}}^{-}}\right]  \right)  ,$
$j\neq j^{{\prime}}.$ In brief, we may say that $Lk\left(  \left[  S_{\alpha
}\right]  \right)  $ contains the infinite bipartite graph with independence
sets $\left \{  \left[  S_{\alpha_{i}^{+}}\right]  ,i=1,2,\ldots \right \}  $ and
$\left \{  \left[  S_{\alpha_{j}^{-}}\right]  ,j=1,2,\ldots \right \}  .$

For each $i,j$ the vertex $\left[  Q_{\alpha_{i}^{+},\alpha_{j}^{-}}\right]  $
is connected with the vertices $\left[  S_{\alpha_{i}^{+}}\right]  ,$ $\left[
S_{\alpha_{j}^{-}}\right]  ,$ $\left[  P_{\alpha,\alpha_{i}^{+}}\right]  ,$
$\left[  P_{\alpha,\alpha_{j}^{-}}\right]  ,$ $\left[  \Sigma_{\alpha_{i}^{+}%
}\right]  $ and $\left[  \Sigma_{\alpha_{j}^{-}}\right]  .$ Thus, there exist
infinitely many vertices of valence $6$ in $Lk\left(  \left[  S_{\alpha
}\right]  \right)  .$ Moreover, it can be checked that all the other vertices
in $Lk\left(  \left[  S_{\alpha}\right]  \right)  $ are of infinite valence.

We now state, for later use, the above three properties for $Lk\left(  \left[
S_{\alpha}\right]  \right)  $ when $S_{\alpha}$ is an annular surface with
separating boundary. We will denote by $K_{n}$ the complete graph on $n$ vertices.

\begin{description}
\item[($sA$-$1$)] $Lk\left(  \left[  S_{\alpha}\right]  \right)  $ contains as
a subgraph the bi-infinite complete bipartite graph,

\item[($sA$-$2$)] There exist infinitely many vertices of valence $6$ in
$Lk\left(  \left[  S_{\alpha}\right]  \right)  $ and all other vertices are of
infinite valence.

\item[($sA$-$3$)] Let $K_{i,j}$ denote the complete graph on the following $6$
vertices: $\left[  S_{\alpha_{i}^{+}}\right]  ,$ $\left[  S_{\alpha_{j}^{-}%
}\right]  ,$ $\left[  P_{\alpha,\alpha_{i}^{+}}\right]  ,$ $\left[
P_{\alpha,\alpha_{j}^{-}}\right]  ,$ $\left[  \Sigma_{\alpha_{i}^{+}}\right]
$ and $\left[  Q_{\alpha_{i}^{+},\alpha_{j}^{-}}\right]  .$ For $i\neq
i^{\prime}$ and $j\neq j^{\prime},$ $K_{i,j}$ and $K_{i^{\prime},j^{\prime}}$
are subgarphs of $Lk\left(  \left[  S_{\alpha}\right]  \right)  $ isomorphic
to $K_{6}$ with no common vertex.

\item[($sA$-$4$)] $Lk\left(  \left[  S_{\alpha}\right]  \right)  $ contains
simple closed cycles of length $3.$
\end{description}

\subsection{Annular surfaces with non-separating boundary}

Let $\alpha$ be a non-separating curve in $\partial M$ which is not
homotopically trivial in $M.$ Let $S_{\alpha}$ be the properly embedded
annular (incompressible) surface whose boundary consists of two (parallel)
copies of $\alpha$ so that $S_{\alpha}$ is isotopic to the annular region in
$\partial M$ bounded by its boundary components. We consider two cases
according to whether or not $\alpha$ intersects all meridians.

\subsubsection{Annular surfaces with non-separating boundary which intersects
all meridians\label{nsa1}}

The incompressible surfaces which give rise to vertices in $Lk\left(  \left[
S_{\alpha}\right]  \right)  $ can be divided into two classes:

\begin{itemize}
\item surfaces $S$ with $\partial S\cap a=\varnothing$ so that $\partial
M\setminus \partial S$ does not contain a separating (for $\partial M$) curve

\item surfaces $S$ with $\partial S\cap a=\varnothing$ so that $\partial
M\setminus \partial S$ contains a separating (for $\partial M$) curve
\end{itemize}

Each surface $S$ in the former class is, necessarily, a pair of pants with all
three boundary components being non-separating, mutually non-isotopic
essential curves with one boundary components of $S$ being isotopic to
$\alpha.$ These surfaces can be enumerated as follows: consider the infinite
collection
\[
\left \{  \left \{  \left[  \delta_{i}\right]  ,\left[  \delta_{j}\right]
\right \}  \bigm \vert i,j=0,1,2,\ldots \right \}
\]
of all distinct (unordered) pairs of isotopy classes of essential curves
$\delta_{i},\delta_{j}$ such that all $\delta_{i},\delta_{j}$ are non-meridian
and non-separating, for each pair $\left \{  \left[  \delta_{i}\right]
,\left[  \delta_{j}\right]  \right \}  $ the curves $\alpha,\delta_{i}%
,\delta_{j}$ are mutually non-isotopic and $\left[  \delta_{i}\right]
\cap \left[  a\right]  =\varnothing,$ $\left[  \delta_{j}\right]  \cap \left[
a\right]  =\varnothing$, $\left[  \delta_{i}\right]  \cap \left[  \delta
_{j}\right]  =\varnothing.$ Each such pair, gives rise to two pairs of pants
$P_{\alpha,\delta_{i},\delta_{j}}^{+},$ $P_{\alpha,\delta_{i},\delta_{j}}^{-}$
with boundary components being isotopic to $\alpha,\delta_{i},\delta_{j}$
respectively. As $\delta_{i},\delta_{j}$ do not bound a disk, both surface are
incompressible and, by Lemma \ref{sigma_torus}, give rise to distinct vertices
$\left[  P_{\alpha,\delta_{i},\delta_{j}}^{+}\right]  ,$ $\left[
P_{\alpha,\delta_{i},\delta_{j}}^{-}\right]  $ in $Lk\left(  \left[
S_{\alpha}\right]  \right)  .$ To see that they are distinct, observe that if
$P_{\alpha,\delta_{i},\delta_{j}}^{+},$ $P_{\alpha,\delta_{i},\delta_{j}}^{-}$
were isotopic, then $M$ would be homeomorphic to $P_{\alpha,\delta_{i}%
,\delta_{j}}^{+}\times \left[  0,1\right]  $ and, thus, we may find a meridian
not intersecting $\alpha.$

We proceed now with the surfaces in the second class, namely, those $S$ for
which $\partial M\setminus \partial S$ contains a separating (for $\partial M$)
curve. The curve $\alpha$ determines a sequence $\left \{  \left[  \beta
_{i}\right]  ,i=1,2,\ldots \right \}  $ consisting of all isotopy classes of
separating curves with the property $\left[  \alpha \right]  \cap \left[
\beta_{i}\right]  =\varnothing$ . This can be done by enumerating the isotopy
classes, say, $\left[  \gamma_{i}\right]  ,i=1,2,\ldots$ of simple closed
curves which intersect $\left[  \alpha \right]  $ at exactly one point and then
taking $\beta_{i}$ to be the commutator $\gamma_{i}a\gamma_{i}^{-1}\alpha
^{-1}.$ Note that the sequence $\left \{  \left[  \beta_{i}\right]
,i=1,2,\ldots \right \}  $ does not contain meridian separating curves since all
meridians intersect $\alpha.$ We will be writing $\left[  S_{\beta_{i}%
}\right]  $ for the corresponding annular vertex which clearly belongs to
$Lk\left(  \left[  S_{\alpha}\right]  \right)  .$ We will complete the full
list of vertices in $Lk\left(  \left[  S_{\alpha}\right]  \right)  $ by
looking at all incompressible surfaces $S$ in $M$ whose boundary does not
intersect $\alpha$ nor $\beta_{i}$ for a fixed $i.$ Obviously, each such
incompressible surface is connected by an edge with the annular vertex
$S_{\beta_{i}}.$ We will then let $i$ vary.

Next we fix a separating non-meridian curve $\beta_{i}$ along with the
corresponding annular surface $\left[  S_{\beta_{i}}\right]  .$

\begin{itemize}
\item As $\beta_{i}$ intersects all meridians, $T_{\beta_{i},+}$ and
$T_{\beta_{i},-}$ are incompressible surfaces (we agree that $T_{\beta_{i},+}$
contains $\alpha$) and they determine vertices $\left[  T_{\beta_{i}%
,+}\right]  ,\left[  T_{\beta_{i},-}\right]  $ in $Lk\left(  \left[
S_{\alpha}\right]  \right)  .$ Recall that these two vertices are not
necessarily distinct (see Remark \ref{r4})$.$

\item Let $\beta_{i}^{j},j=1,2,\ldots$ be the infinite sequence of all
pair-wise non-isotopic essential simple closed curves in $T_{\beta_{i},-}$
which intersect pair-wise. The corresponding annular surfaces $S_{\beta
_{i}^{j}},$ $j=1,2,\ldots$ give rise to distinct vertices $\left[
S_{\beta_{i}^{j}}\right]  $ for all $j.$

\item Since $\alpha$ intersects all meridians, the surface $\Sigma_{\alpha}$
is incompressible and gives rise to a vertex $\left[  \Sigma_{\alpha}\right]
$ in $Lk\left(  \left[  S_{\alpha}\right]  \right)  .$

\item Each surface $\Sigma_{\beta_{i}^{j}}$ may or may not be incompressible
depending on whether, or not, $\beta_{i}^{j}$ intersects all meridians.
However, by Proposition \ref{algebraic}, for each $i,$ there exist infinitely
many $j$'s such that $\beta_{i}^{j}$ intersects all meridians and, hence, we
obtain distinct vertices $\left[  \Sigma_{\beta_{i}^{j}}\right]  $ in
$Lk\left(  \left[  S_{\alpha}\right]  \right)  $ for infinitely many $j$'s.

\item The surface $P_{\beta_{i},\alpha}$ is incompressible and it gives rise
to a vertex $\left[  P_{\beta_{i},\alpha}\right]  $ in $Lk\left(  \left[
S_{\alpha}\right]  \right)  .$

\item For each $j=1,2,\ldots$ the surface $P_{\beta_{i},\beta_{i}^{j}}$ is
incompressible because $\beta_{i}$ is separating and non meridian, hence, by
Proposition \ref{all_meridians}, intersects all meridians. Thus, we obtain
distinct vertices $\left[  P_{\beta_{i},\beta_{i}^{j}}\right]  $ in $Lk\left(
\left[  S_{\alpha}\right]  \right)  $ for all $j.$

\item For all $j=1,2,\ldots$ the surfaces $Q_{\alpha,\beta_{i}^{j}}$ are
incompressible. These surfaces give rise to distinct vertices $\left[
Q_{\alpha,\beta_{i}^{j}}\right]  ,j=1,2,\ldots.$
\end{itemize}

As in the previous section, we may verify that the above list is a complete
list of vertices in $Lk\left(  \left[  S_{\alpha}\right]  \right)  $ which are
connected by an edge with the fixed vertex $\left[  S_{\beta_{i}}\right]  .$

It is clear that each vertex $\left[  P_{\beta_{i},\beta_{i}^{j}}\right]  $ is
connected by an edge with the vertices $\left[  S_{\beta_{i}}\right]  ,$
$\left[  S_{\beta_{i}^{j}}\right]  ,$ $\left[  \Sigma_{\alpha}\right]  ,$
$\left[  \Sigma_{\beta_{i}^{j}}\right]  ,$ $\left[  P_{\beta_{i},\alpha
}\right]  ,$ $\left[  Q_{\alpha,\beta_{i}^{j}}\right]  ,$ $\left[
T_{\beta_{i},+}\right]  ,$ $\left[  T_{\beta_{i},-}\right]  .$ Thus,
infinitely many vertices have valence $8$ or, in the case $\left[
T_{\alpha,+,i}\right]  =\left[  T_{\alpha,+,i}\right]  ,$ valence $7.$ All
other vertices in $Lk\left(  \left[  S_{\alpha}\right]  \right)  $ have
infinite valence.

We conclude this section by stating the above properties of $Lk\left(  \left[
S_{\alpha}\right]  \right)  $ when $S_{\alpha}$ is an annular surface with non
separating boundary ($nsA$)$.$ These properties will be used later in Section
\ref{invariance_proof} to prove Proposition \ref{non_iso}.

\begin{description}
\item[($nsA$-$1$)] $Lk\left(  \left[  S_{\alpha}\right]  \right)  $ contains
infinitely many vertices of valence $7$ or $8.$

\item[($nsA$-$2$)] Infinitely many vertices in $Lk\left(  \left[  S_{\alpha
}\right]  \right)  $ have infinite valence.

\item[($nsA$-$3$)] $Lk\left(  \left[  S_{\alpha}\right]  \right)  $ contains
simple closed cycles of length $3.$
\end{description}

\subsubsection{Annular surfaces with non-separating boundary which does not
intersect all meridians}

In this case the list of vertices is similar to the previous subsection, the
main difference being that $\Sigma_{\alpha}$ is not incompressible, hence,
$\left[  \Sigma_{\alpha}\right]  $ does not exist in $Lk\left(  \left[
S_{\alpha}\right]  \right)  .$ It is clear that properties ($nsA$-$2$),
($nsA$-$3$) are still valid and we proceed to show that ($nsA$-$1$) holds in
this case as well.

There exist infinitely many separating non-meridian curves $\left \{  \left[
\beta_{i}\right]  ,i=1,2,\ldots \right \}  $ which do not intersect $\alpha.$
For each $i,$ using Proposition \ref{algebraic}, we obtain non-separating
curves $\beta_{i}^{j}$ with $\beta_{i}^{j}\cap \beta_{i}=\varnothing$ such that
$\beta_{i}^{j}$ intersects all meridians for infinitely many $j$'s As in the
previous subsection, each vertex $\left[  P_{\beta_{i},\beta_{i}^{j}}\right]
$ is connected by an edge with the vertices $\left[  S_{\beta_{i}}\right]  ,$
$\left[  S_{\beta_{i}^{j}}\right]  ,$ $\left[  \Sigma_{\beta_{i}^{j}}\right]
,$ $\left[  P_{\beta_{i},\alpha}\right]  ,$ $\left[  Q_{\alpha,\beta_{i}^{j}%
}\right]  ,$ $\left[  T_{\beta_{i},+}\right]  ,$ $\left[  T_{\beta_{i}%
,-}\right]  $ and its valence is $7$ provided that $\left[  T_{\beta_{i}%
,+}\right]  \neq \left[  T_{\beta_{i},-}\right]  .$

Thus, in order to establish property ($nsA$-$1$) for this case it suffices to
show that for infinitely many $i$'s the surfaces $T_{\beta_{i},+},$
$T_{\beta_{i},-}$ are not isotopic. To do this, we first claim that if
$\alpha$ is not a generator of $\pi_{1}\left(  M\right)  $ then $\left[
T_{\beta_{i},+}\right]  \neq \left[  T_{\beta_{i},-}\right]  $. Assume that for
a separating curve $\beta_{i_{0}},$ the surfaces $T_{\beta_{i_{0}},+}$ and
$T_{\beta_{i_{0}},-}$ are isotopic (we agree that $T_{\beta_{i_{0}},+}$
contains $\alpha$). Choose a curve $\gamma$ such that the commutator $\left[
\alpha,\gamma \right]  =\beta_{i_{0}}.$ Equivalently, $\alpha$ and $\gamma$
intersect at $1$ point and $\alpha,\gamma$ generate $\pi_{1}\left(
T_{\beta_{i},+}\right)  .$ As $T_{\beta_{i_{0}},+}$ and $T_{\beta_{i_{0}},-}$
are isotopic, $M$ is homeomorphic to $T_{\beta_{i},+}\times \left[  0,1\right]
$ and, hence, the generators of $\pi_{1}\left(  T_{\beta_{i},+}\right)  $
generate $\pi_{1}\left(  M\right)  ,$ a contradiction.

Thus, without loss of generality, we assume that $\alpha$ is a generator of
$\pi_{1}\left(  M\right)  .$ Find a curve $x$ on $\partial M\ $such that
$x\cap \alpha=\varnothing$ and the corresponding element $x,\alpha$ generate
$\pi_{1}\left(  M\right)  .$ The curves $x^{2}$ and $x\alpha$ intersect at one
point and the commutator $\beta_{2}=\left[  x^{2},x\alpha \right]  $ is a
separating curve in $\partial M$ and does not bound a disk in $M.$ We claim
that the surfaces $T_{\beta_{2},+},$ $T_{\beta_{2},-}$ corresponding to the
separating curve $\beta_{2}$ are not isotopic. For, if $T_{\beta_{2},+}$ is
isotopic to $T_{\beta_{2},-}$ then $M$ would be homeomorphic to $T_{\beta
_{2},-}\times \left[  0,1\right]  .$ In particular, any generator of $\pi
_{1}\left(  T_{\beta_{2},-}\right)  $ would be generator for $\pi_{1}\left(
M\right)  .$ As $x^{2}$ and $x\alpha$ are generators for $T_{\beta_{2},-}$ we
would have that $x^{2}$ is a generator for $\pi_{1}\left(  M\right)  .$ This
is a contradiction since $x^{2}$ is not a generator for the free group of rank
$2$ when $x$ is. In a similar manner and using the fact that $x^{i},i\geq2$ is
not a generator for the free group of rank $2$ when $x$ is, we construct
infinitely many curves $\beta_{i}=\left[  x^{i},x^{i-1}\alpha \right]  $ such
that the corresponding surfaces $T_{\beta_{i},+},$ $T_{\beta_{i},-}$ are not
isotopic. This completes the proof of existence of infinitely many vertices of
valence $7.$

\subsection{Surfaces of type $\left(  T\right)  $}

Let $\alpha$ be a separating curve in $\partial M$ which is not homotopically
trivial in $M.$ Denote by $T_{\alpha,+},T_{\alpha,-}$ the closures of the
components of $\partial M\setminus S_{\alpha}$ (each being a torus with one
boundary component). We will study the link of $\left[  T_{\alpha,+}\right]
.$ Note that $T_{\alpha,-}$ may or may not be isotopic to $T_{\alpha,+}.$

As before, denote by $\alpha_{i}^{+},i=1,2,\ldots$ (resp. $\alpha_{j}%
^{-},j=1,2,\ldots$) the infinite sequence of all pair-wise non-isotopic
essential simple closed curves in $T_{\alpha,+}$ (resp. $T_{\alpha,-}$) which
intersect pair-wise.

\begin{itemize}
\item The corresponding annular surfaces $S_{\alpha_{i}^{+}},$ $i=1,2,\ldots$
and $S_{\alpha_{j}^{-}}$, $j=1,2,\ldots$ give rise to distinct vertices
$\left[  S_{\alpha_{i}^{+}}\right]  ,\left[  S_{\alpha_{j}^{-}}\right]  $
which all belong to $Lk\left(  \left[  T_{\alpha,+}\right]  \right)  .$

\item The surface $\Sigma_{\alpha_{i}^{+}}$ (resp. $\Sigma_{\alpha_{j}^{-}}$)
may or may not be incompressible depending on whether, or not, $\alpha_{i}%
^{+}$ (resp. $\alpha_{j}^{-}$) intersects all meridians. By Proposition
\ref{algebraic}, there exist infinitely many $i$'s ( resp. $j$'s) such that
$\alpha_{i}^{+}$ (resp. $\alpha_{j}^{-}$) intersects all meridians. By Lemma
\ref{sigma_torus}, $\left[  \Sigma_{\alpha_{j}^{-}}\right]  $ belongs to
$Lk\left(  \left[  T_{\alpha,+}\right]  \right)  $ for infinitely many $j$'s
whereas none of the surfaces $\Sigma_{\alpha,\alpha_{i}^{+}}$ gives rise to a
vertex in $Lk\left(  \left[  T_{\alpha,+}\right]  \right)  .$

\item The surface $P_{\alpha,\alpha_{i}^{+}}$ (resp. $P_{\alpha,\alpha_{j}%
^{-}}$) is incompressible for all $i$ (resp. $j$). These surfaces give rise to
distinct vertices $\left[  P_{\alpha,\alpha_{i}^{+}}\right]  ,i=1,2,\ldots$
and $\left[  P_{\alpha,\alpha_{j}^{-}}\right]  ,j=1,2,\ldots$ in $Lk\left(
\left[  T_{\alpha,+}\right]  \right)  .$
\end{itemize}

These are all the vertices in $Lk\left(  \left[  T_{\alpha,+}\right]  \right)
.$ If $T_{\alpha,-}$ is isotopic to $T_{\alpha,+}$ then the vertex $\left[
T_{\alpha,-}\right]  $ is not present and all other vertices mentioned above
exist. Moreover, $\left[  \Sigma_{\alpha_{j}^{-}}\right]  $ belongs to
$Lk\left(  \left[  T_{\alpha,+}\right]  \right)  $ for infinitely many $i$'s.
We conclude this section by stating two easily checked properties for
$Lk\left(  \left[  T\right]  \right)  $ when $T$ is a surface of type $\left(
T\right)  .$

\begin{description}
\item[($T$-$1$)] All vertices in $Lk\left(  \left[  T\right]  \right)  $ are
of infinite valence.

\item[($T$-$2$)] $Lk\left(  \left[  T\right]  \right)  $ contains simple
closed cycles of length $3.$
\end{description}

\subsection{Surfaces of type $\left(  \Sigma \right)  $}

Let $\alpha$ be a non-separating curve in $\partial M$ which is not
homotopically trivial in $M$ such that $\alpha$ intersects all meridians. Such
a curve defines a genus $1$ incompressible surface $\Sigma_{\alpha}$ with two
boundary components homotopic to $\alpha.$ As before, the curve $\alpha$
determines an infinite sequence $\left \{  \left[  \beta_{i}\right]
,i=1,2,\ldots \right \}  $ of isotopy classes of separating curves with the
property $\left[  \alpha \right]  \cap \left[  \beta_{i}\right]  =\varnothing$ .
As $\Sigma_{\alpha}$ is incompressible, none of the $\beta_{i}$'s bounds a
disk. The incompressible surfaces which give rise to vertices in $Lk\left(
\left[  \Sigma_{\alpha}\right]  \right)  $ can be divided into two classes:

\begin{itemize}
\item surfaces $S$ with $\partial S\cap a=\varnothing$ so that $\partial
M\setminus \partial S$ does not contain a separating (for $\partial M$) curve

\item surfaces $S$ with $\partial S\cap a=\varnothing$ so that $\partial
M\setminus \partial S$ contains a separating (for $\partial M$) curve
\end{itemize}

Note that a surface $S$ in the former class is, necessarily, a pair of pants
with all three boundary components being non separating curves (in $\partial
M$).

Taking into account the above two classes of surfaces, we will compose the
full list of vertices of $Lk\left(  \left[  \Sigma_{\alpha}\right]  \right)  $
by looking at

\begin{itemize}
\item[(A)] all incompressible surfaces $S$ which are pair of pants with all
three boundary components being non-separating, mutually non-isotopic
essential curves with one boundary components of $S$ being isotopic to
$\alpha,$ and

\item[(B)] all incompressible surfaces $S$ in $M$ whose boundary does not
intersect $\alpha$ nor $\beta_{i}$ for a fixed $i.$ Obviously, each such
incompressible surface is connected by an edge with the annular surface
$S_{\beta_{i}}$ (recall, $\beta_{i}$ is not a meridian curve). We will then
let $i$ vary.
\end{itemize}

The surfaces in the former class which give rise to vertices in $Lk\left(
\left[  \Sigma_{\alpha}\right]  \right)  $ have been analyzed in the previous
section: they are pairs of pants $P_{\alpha,\delta_{i},\delta_{j}}^{+},$
$P_{\alpha,\delta_{i},\delta_{j}}^{-}$ with boundary components being isotopic
to $\alpha,\delta_{i},\delta_{j}$ respectively which can be enumerated by the
infinite collection
\[
\left \{  \left \{  \left[  \delta_{i}\right]  ,\left[  \delta_{j}\right]
\right \}  \bigm \vert i,j=0,1,2,\ldots \right \}
\]
of all distinct (unordered) pairs of isotopy classes of essential curves
$\delta_{i},\delta_{j}$ such that: $\delta_{i},\delta_{j}$ are non-meridian
and non-separating curves, the curves $\alpha,\delta_{i},\delta_{j}$ are
mutually non-isotopic and $\left[  \delta_{i}\right]  \cap \left[  a\right]
=\varnothing,$ $\left[  \delta_{j}\right]  \cap \left[  a\right]
=\varnothing,$ $\left[  \delta_{j}\right]  \cap \left[  \delta_{i}\right]
=\varnothing.$ As $\delta_{i},\delta_{j}$ do not bound a disk, each such
surface is incompressible and, by Lemma \ref{sigma_torus}, gives rise to
(necessarily distinct, as explained in Subsection \ref{nsa1}) vertices
$\left[  P_{\alpha,\delta_{i},\delta_{j}}^{+}\right]  ,$ $\left[
P_{\alpha,\delta_{i},\delta_{j}}^{-}\right]  $ in $Lk\left(  \left[
\Sigma_{\alpha}\right]  \right)  .$

For the surfaces in the second class which give rise to vertices in $Lk\left(
\left[  \Sigma_{\alpha}\right]  \right)  ,$ we will fix a separating curve
$\beta_{i}$ and look at all incompressible surfaces in $M$ which, up to
isotopy, do not intersect $\beta_{i},$ for a fixed $i.$ We shall then let $i$
vary in order to complete the list of vertices of $Lk\left(  \left[
\Sigma_{\alpha}\right]  \right)  .$

Fix a separating curve $\beta_{i}$ and the corresponding annular surface
$\left[  S_{\beta_{i}}\right]  .$

\begin{itemize}
\item Clearly, $\left[  S_{\alpha}\right]  $ as well as all $\left[
S_{\beta_{i}}\right]  $ belong to $Lk\left(  \left[  \Sigma_{\alpha}\right]
\right)  .$

\item Let $T_{\beta_{i},-}$ be the surface not containing $\alpha.$ As
$T_{\beta_{i},-}$ is incompressible, $\left[  T_{\beta_{i},-}\right]  $
belongs to $Lk\left(  \left[  \Sigma_{\alpha}\right]  \right)  $ by Lemma
\ref{sigma_torus}. Observe that $T_{\beta_{i},+}$ is not, up to isotopy,
disjoint from $\Sigma_{\alpha},$ hence, $\left[  T_{\beta_{i},+}\right]  $
does not exist in $Lk\left(  \left[  \Sigma_{\alpha}\right]  \right)  .$ As
$\alpha$ intersects all meridians, it can be shown that $T_{\beta_{i},-}$ is
never isotopic to $T_{\beta_{i},+}.$

\item Let $\beta_{i}^{j},j=1,2,\ldots$ be the infinite sequence of all
essential simple closed curves in $T_{\beta_{i},-}$ which, pair-wise,
intersect and are non-isotopic. The corresponding annular surfaces
$S_{\beta_{i}^{j}},$ $j=1,2,\ldots$ give rise to distinct vertices $\left[
S_{\beta_{i}^{j}}\right]  $ in $Lk\left(  \left[  \Sigma_{\alpha}\right]
\right)  .$

\item The surface $P_{\beta_{i},\alpha}$ is incompressible and it gives rise
to a vertex $\left[  P_{\beta_{i},\alpha}\right]  $ in $Lk\left(  \left[
S_{\alpha}\right]  \right)  .$

\item For each $j=1,2,\ldots$ the surface $P_{\beta_{i},\beta_{i}^{j}}$ is
incompressible because $\beta_{i}$ is separating and non meridian, hence, by
Proposition \ref{all_meridians}, intersects all meridians. Thus, we obtain
distinct vertices $\left[  P_{\beta_{i},\beta_{i}^{j}}\right]  $ in $Lk\left(
\left[  S_{\alpha}\right]  \right)  $ for all $j.$

\item For each $j=1,2,\ldots$ the surface $Q_{\alpha,\beta_{i}^{j}}$ is
incompressible. These surfaces give rise to distinct vertices $\left[
Q_{\alpha,\beta_{i}^{j}}\right]  ,j=1,2,\ldots$ in $Lk\left(  \left[
\Sigma_{\alpha}\right]  \right)  .$
\end{itemize}

These are all the vertices in $Lk\left(  \left[  \Sigma_{\alpha}\right]
\right)  .$

For all $i,j$ the vertex $\left[  P_{\alpha,\delta_{i},\delta_{j}}^{+}\right]
$ is connected by an edge with the following vertices: $\left[  S_{\alpha
}\right]  ,$ $\left[  Q_{\alpha,\delta_{i}}\right]  ,$ $\left[  Q_{\alpha
,\delta_{j}}\right]  ,$ $\left[  S_{\delta_{i}}\right]  ,$ $\left[
S_{\delta_{j}}\right]  $ and $\left[  P_{\alpha,\delta_{i},\delta_{j}}%
^{-}\right]  .$ Recall that $\left[  P_{\alpha,\delta_{i},\delta_{j}}%
^{+}\right]  \neq \left[  P_{\alpha,\delta_{i},\delta_{j}}^{-}\right]  .$ Thus,
each vertex $\left[  P_{\alpha,\delta_{i},\delta_{j}}^{+}\right]  $ has
valence $6.$

Similarly, we check that for all $i,j$ the vertex $\left[  P_{\beta_{i}%
,\alpha_{j}^{i}}\right]  $ is connected by the following vertices: $\left[
S_{\alpha}\right]  ,$ $\left[  S_{\beta_{i}}\right]  ,$ $\left[  T_{\beta
_{i},-}\right]  ,$ $\left[  S_{\alpha_{j}^{i}}\right]  ,$ $\left[
Q_{\alpha,\alpha_{j}^{i}}\right]  $ and $\left[  P_{\beta_{i},\alpha}\right]
.$ Thus, all vertices $\left[  P_{\beta_{i},\alpha_{j}^{i}}\right]  $ have
valence $6.$ All the remaining vertices in $Lk\left(  \left[  \Sigma_{\alpha
}\right]  \right)  $ have infinite valence. For example, for fixed $i,j$ the
surface $Q_{\alpha,\alpha_{j}^{i}}$ contains infinitely many separating curves
hence, the vertex $\left[  Q_{\alpha,\alpha_{j}^{i}}\right]  $ is connected by
$\left[  S_{\beta_{i}}\right]  $ for infinitely many $i$'s.

Moreover, it is easy to check that any subgraph of $Lk\left(  \left[
\Sigma_{\alpha}\right]  \right)  $ isomorphic to $K_{6}$ either contains a
vertex of the form $\left[  P_{\alpha,\delta_{i},\delta_{j}}^{+}\right]  $ for
some $i,j$ or, a vertex of the form $\left[  P_{\beta_{i},\alpha_{j}^{i}%
}\right]  $ for some $i,j.$ In the former case, the remaining five vertices
are $\left[  S_{\alpha}\right]  ,$ $\left[  S_{\delta_{i}}\right]  ,$ $\left[
S_{\delta_{j}}\right]  ,$ $\left[  P_{\alpha,\delta_{i},\delta_{j}}%
^{-}\right]  $ and one of $\left[  Q_{\alpha,\delta_{i}}\right]  ,$ $\left[
Q_{\alpha,\delta_{j}}\right]  .$ In the latter case the remaining five
vertices are $\left[  S_{\alpha}\right]  ,$ $\left[  S_{\beta_{i}}\right]  ,$
$\left[  S_{\alpha_{j}^{i}}\right]  ,$ $\left[  P_{\beta_{i},\alpha}\right]  $
and one of $\left[  T_{\beta_{i},-}\right]  ,$ $\left[  Q_{\alpha,\alpha
_{j}^{i}}\right]  .$ In both cases, any subgraph of $Lk\left(  \left[
\Sigma_{\alpha}\right]  \right)  $ isomorphic to $K_{6}$ must contain $\left[
S_{\alpha}\right]  .$

We conclude this section by stating two properties for $Lk\left(  \left[
\Sigma \right]  \right)  $ when $\Sigma$ is a surface of type $\left(
\Sigma \right)  .$

\begin{description}
\item[($\Sigma$-$1$)] There exist infinitely many vertices of valence $6$ in
$Lk\left(  \left[  \Sigma \right]  \right)  .$

\item[($\Sigma$-$2$)] Any subgraph of $Lk\left(  \left[  \Sigma_{\alpha
}\right]  \right)  $ isomorphic to $K_{6}$ must contain $\left[  S_{\alpha
}\right]  .$

\item[($\Sigma$-$3$)] $Lk\left(  \left[  \Sigma \right]  \right)  $ contains
simple closed cycles of length $3.$
\end{description}

\subsection{Surfaces of type $\left(  P\right)  $}

If $P$ is a pair of pants properly embedded and incompressible in $M$ with
boundary components $\alpha,\beta,\gamma$ each being non-separating and
mutually non-isotopic. As $P$ is incompressible, none of the curves
$\alpha,\beta,\gamma$ bounds a disk. Thus, the annular surfaces $S_{\alpha
},S_{\beta},S_{\gamma}$ are incompressible and the corresponding vertices
$\left[  S_{\alpha}\right]  ,\left[  S_{\beta}\right]  ,\left[  S_{\gamma
}\right]  $ belong to $Lk\left(  \left[  P\right]  \right)  .$ Each surfaces
$\Sigma_{\alpha},$ $\Sigma_{\beta},$ $\Sigma_{\gamma}$ may or may not be
incompressible depending on whether, or not, each of the curves $\alpha,$
$\beta,$ $\gamma$ intersects all meridians. Thus, the vertices $\left[
\Sigma_{\alpha}\right]  ,\left[  \Sigma_{\beta}\right]  ,\left[
\Sigma_{\gamma}\right]  $ may or may not exist in $Lk\left(  \left[  P\right]
\right)  .$ Similarly, the vertices $\left[  Q_{\alpha,\beta}\right]  ,\left[
Q_{\beta,\gamma}\right]  ,\left[  Q_{a,\gamma}\right]  $ may or may not exist
in $Lk\left(  \left[  P\right]  \right)  .$ In any case, $Lk\left(  \left[
P\right]  \right)  $ contains finitely many vertices.

We proceed now with the case $\partial P$ contains a separating curve. Let
$\alpha$ be a non-separating curve in $\partial M$ and $\beta$ a separating
curve in $\partial M$ with $\alpha \cap \beta=\varnothing,$ both not
homotopically trivial in $M.$ Let $P_{\beta,\alpha}$ be the pair of pants with
boundary components $\alpha,\alpha$ and $\beta.$

\begin{itemize}
\item Clearly, $\left[  S_{\alpha}\right]  ,$ $\left[  S_{\beta}\right]  $
belong to $Lk\left(  \left[  P_{\beta,\alpha}\right]  \right)  .$

\item Let $T_{\beta,+}$ (resp. $T_{\beta,-}$) be the surface which contains
(resp. does not contain) $\alpha.$ Both $T_{\beta,-},$ $T_{\beta,+}$ are
incompressible and give rise to distinct vertices in $Lk\left(  \left[
P_{\alpha,\beta}\right]  \right)  $ unless $T_{\beta,+}$ is isotopic to
$T_{\beta,-},$ in which case, $\left[  T_{\beta,+}\right]  =\left[
T_{\beta,-}\right]  .$

\item Denote by $\alpha_{i},i=1,2,\ldots$ the infinite sequence of pair-wise
non-isotopic essential simple closed curves in $T_{\beta,-}$ which intersect
pair-wise. The corresponding annular surfaces $S_{\alpha_{i}},$ $i=1,2,\ldots$
give rise to distinct vertices $\left[  S_{\alpha_{i}}\right]  $ in $Lk\left(
\left[  P_{\beta,\alpha}\right]  \right)  .$

\item The surfaces $P_{\beta,\alpha_{i}}$ are incompressible for all $i$ and
give rise to distinct vertices $\left[  P_{\beta,\alpha_{i}}\right]
,i=1,2,\ldots$ in $Lk\left(  \left[  P_{\beta,\alpha}\right]  \right)  .$

\item The surface $\Sigma_{\alpha}$ may or may not be incompressible, hence,
$\left[  \Sigma_{\alpha}\right]  $ may or may not exist in $Lk\left(  \left[
P_{\beta,\alpha}\right]  \right)  .$

\item For each $i,$ the surface $\Sigma_{\alpha_{i}}$ may or may not be
incompressible depending on whether, or not, $\alpha_{i}$ intersects all
meridians. By Lemma \ref{algebraic}, there exist infinitely many $i$'s such
that $\alpha_{i}$ intersects all meridians and, hence, we obtain distinct
vertices $\left[  \Sigma_{\alpha_{i}}\right]  $ in $Lk\left(  \left[
P_{\beta,\alpha}\right]  \right)  $ for infinitely many $i$'s.

\item For each $i=1,2,\ldots$ the surface $Q_{\alpha,\alpha_{i}}$ may or may
not be incompressible depending on whether, or not, $\alpha_{i}$ intersects
all meridians. By Lemma \ref{algebraic}, there exist infinitely many $i$'s
such that $\alpha_{i}$ intersects all meridians and, hence, we obtain distinct
vertices $\left[  Q_{\alpha,\alpha_{i}}\right]  $ in $Lk\left(  \left[
P_{\beta,\alpha}\right]  \right)  $ for infinitely many $i$'s.
\end{itemize}

These are all the vertices in $Lk\left(  \left[  P_{\alpha,\beta}\right]
\right)  $.

We conclude this section by stating properties for $Lk\left(  \left[
P\right]  \right)  $ when $P$ is a surface of type $\left(  P\right)  $ when
$\partial P$ contains a separating curve.

\begin{description}
\item[($P$-$1$)] $Lk\left(  \left[  P\right]  \right)  $ contains finitely
many vertices of infinite valence, namely, the vertices $\left[  S_{\alpha
}\right]  ,$ $\left[  S_{\beta}\right]  ,$ $\left[  \Sigma_{\alpha}\right]  ,$
$\left[  T_{\beta,+}\right]  $ and/or $\left[  T_{\beta,-}\right]  .$

\item[($P$-$2$)] For infinitely many $i$'s, $Lk\left(  \left[  P\right]
\right)  $ contains $K_{6}$ as a subgraph consisting of the following
vertices: $\left[  S_{\alpha}\right]  ,$ $\left[  S_{\beta}\right]  ,$
$\left[  S_{\alpha_{i}}\right]  ,$ $\left[  P_{\beta,\alpha_{i}}\right]  ,$
$\left[  \Sigma_{\alpha_{i}}\right]  ,$ $\left[  Q_{\alpha,\alpha_{i}}\right]
.$
\end{description}

Observe that the presence or absence of $\left[  \Sigma_{\alpha}\right]  $ in
$Lk\left(  \left[  P_{\beta,\alpha}\right]  \right)  $ does not affect the
above two properties. If $P$ is a pair of pants with all three boundary
components being non-separating, then

\begin{description}
\item[($nsP$)] $Lk\left(  \left[  P\right]  \right)  $ contains finitely many vertices.
\end{description}

\subsection{Surfaces of type $\left(  Q\right)  $}

Let $\alpha,\beta$ be two non-isotopic, non-separating curves in $\partial M$
both not homotopically trivial in $M$ with $\alpha \cap \beta=\varnothing.$ Let
$Q_{\alpha,\beta}$ be the sphere with 4 holes with boundary components
isotopic to $\alpha,\alpha,\beta,$ and $\beta.$ As before, the curves
$\alpha,\beta$ determine

\begin{itemize}
\item the infinite sequence $\left \{  \left[  \gamma_{i}\right]
,i=1,2,\ldots \right \}  $ of isotopy classes of separating curves in $\partial
M$ each having the property $\alpha \cap \gamma_{i}=\varnothing$ and $\beta
\cap \gamma_{i}=\varnothing.$

\item the infinite sequence $\left \{  \left[  \delta_{i}\right]
,i=1,2,\ldots \right \}  $ of isotopy classes of non-separating curves each
having the property $\alpha \cap \delta_{i}=\varnothing$ and $\beta \cap
\delta_{i}=\varnothing.$
\end{itemize}

None of the curves in these classes bounds a disk since $Q_{\alpha,\beta}$ is
incompressible. We will compose the list of vertices of $Lk\left(  \left[
Q_{\alpha,\beta}\right]  \right)  $ by looking at

\begin{itemize}
\item[(A)] all incompressible surfaces in $M$ which, up to isotopy, do not
intersect $\gamma_{i},$ for a fixed $i.$

\item[(B)] all incompressible surfaces in $M$ which, up to isotopy, do not
intersect $\delta_{j},$ for a fixed $j.$
\end{itemize}

We shall then let $i,j$ vary in order to get a complete list of vertices of
$Lk\left(  \left[  Q_{\alpha,\beta}\right]  \right)  .$

Fix a separating curve $\gamma_{i}$ and the corresponding annular surface
$\left[  S_{\gamma_{i}}\right]  .$

\begin{itemize}
\item Apparently, all vertices $\left[  S_{\gamma_{i}}\right]  $ belong to
$Lk\left(  \left[  Q_{\alpha,\beta}\right]  \right)  $ as well as the vertices
$\left[  S_{\alpha}\right]  $ and $\left[  S_{\beta}\right]  .$

\item The surfaces $P_{\gamma_{i},\alpha},$ $P_{\gamma_{i},\beta}$ are
incompressible for all $i$ and give rise to distinct vertices $\left[
P_{\gamma_{i},\alpha}\right]  $ and $\left[  P_{\gamma_{i},\beta}\right]
,i=1,2,\ldots$ in $Lk\left(  \left[  Q_{\alpha,\beta}\right]  \right)  .$

\item The surface $\Sigma_{\alpha}$ (resp. $\Sigma_{\beta}$) may or may not be
incompressible depending on whether, or not, $\alpha$ (resp. $\beta$)
intersects all meridians. Hence, each of $\left[  \Sigma_{\alpha}\right]  $
and $\left[  \Sigma_{\beta}\right]  $ may or may not exist in $Lk\left(
\left[  Q_{\alpha,\beta}\right]  \right)  .$ Observe that if both vertices
$\left[  \Sigma_{\alpha}\right]  ,$ $\left[  \Sigma_{\beta}\right]  $ exist in
$Lk\left(  \left[  Q_{\alpha,\beta}\right]  \right)  $ then they are not
connected by an edge.
\end{itemize}

These are all the vertices in $Lk\left(  \left[  Q_{\alpha,\beta}\right]
\right)  $ which correspond to the class (A) mentioned above.

Now fix a non-separating curve $\delta_{j}$ and the corresponding annular
surface $\left[  S_{\delta_{j}}\right]  .$ All vertices $\left[  S_{\delta
_{j}}\right]  $ belong to $Lk\left(  \left[  Q_{\alpha,\beta}\right]  \right)
.$ Moreover, each curve $\delta_{j}$ gives rise to two pairs of pants
$P_{\alpha,\beta,\delta_{j}}^{+},$ $P_{\alpha,\beta,\delta_{j}}^{-}$ with
boundary components being isotopic to $\alpha,\beta,\delta_{j}$ respectively.
Note that, as every meridian intersects $\alpha \cup \beta,$ $P_{\alpha
,\beta,\delta_{j}}^{+}$ cannot be isotopic to $P_{\alpha,\beta,\delta_{j}}%
^{-}$ (see Subsection \ref{nsa1}). Both surfaces are incompressible for all
$j$ and, by Lemma \ref{sigma_torus}, give rise to distinct vertices $\left[
P_{\alpha,\beta,\delta_{j}}^{+}\right]  ,$ $\left[  P_{\alpha,\beta,\delta
_{j}}^{-}\right]  $ in $Lk\left(  \left[  Q_{\alpha,\beta}\right]  \right)  .$

We conclude this section by stating several properties for $Lk\left(  \left[
Q\right]  \right)  $ when $Q$ is a surface of type $\left(  Q\right)  .$

\begin{description}
\item[($Q$-$1$)] $Lk\left(  \left[  Q\right]  \right)  $ contains finitely
many vertices of infinite valence, namely, the vertices $\left[  S_{\alpha
}\right]  ,$ $\left[  S_{\beta}\right]  ,$ $\left[  \Sigma_{\alpha}\right]  $
and/or $\left[  \Sigma_{\beta}\right]  .$

\item[($Q$-$2.1$)] If neither $\left[  \Sigma_{\alpha}\right]  $ nor $\left[
\Sigma_{\beta}\right]  $ exist in $Lk\left(  \left[  Q\right]  \right)  $
then, $Lk\left(  \left[  Q\right]  \right)  $ contains a subgraph isomorphic
to $K_{5}$ consisting of the following vertices: $\left[  S_{\alpha}\right]
,$ $\left[  S_{\beta}\right]  ,$ $\left[  S_{\gamma_{i}}\right]  ,$ $\left[
P_{\gamma_{i},\alpha}\right]  $ and $\left[  P_{\gamma_{i},\beta}\right]  $
for all $i.$ Moreover, it does not contain a subgraph isomorphic to $K_{6}.$

\item[($Q$-$2.2$)] If exactly one of the vertices $\left[  \Sigma_{\alpha
}\right]  $ and $\left[  \Sigma_{\beta}\right]  $ exists in $Lk\left(  \left[
Q\right]  \right)  ,$ say $\left[  \Sigma_{\alpha}\right]  ,$ then any
subgraph of $Lk\left(  \left[  Q\right]  \right)  $ which is isomorphic to
$K_{6}$ contains the vertices $\left[  S_{\alpha}\right]  ,$ $\left[
S_{\beta}\right]  ,$ $\left[  \Sigma_{\alpha}\right]  $ (the rest $3$ vertices
can be either $\left[  P_{\alpha,\beta,\delta_{j}}^{+}\right]  ,$ $\left[
P_{\alpha,\beta,\delta_{j}}^{-}\right]  ,$ $\left[  S_{\delta_{j}}\right]  $
or, $\left[  S_{\gamma_{i}}\right]  ,$ $\left[  P_{\gamma_{i},\alpha}\right]
,$ $\left[  P_{\gamma_{i},\beta}\right]  $).

\item[($Q$-$2.3$)] If both $\left[  \Sigma_{\alpha}\right]  $ and $\left[
\Sigma_{\beta}\right]  $ exist in $Lk\left(  \left[  Q\right]  \right)  $ then
any subgraph of $Lk\left(  \left[  Q\right]  \right)  $ isomorphic to $K_{6}$
contains the vertices $\left[  S_{\alpha}\right]  ,$ $\left[  S_{\beta
}\right]  $ and exactly one of the vertices $\left[  \Sigma_{\alpha}\right]  $
and $\left[  \Sigma_{\beta}\right]  .$
\end{description}

\subsection{Proof of vertex invariance\label{invariance_proof}}

We are now in position to prove Proposition \ref{non_iso}. Let $D$ be a
separating meridian, $D^{^{\prime}}$ a non-separating meridian, $S_{a}$ an
annular surface with $\alpha$ separating, $S_{a^{\prime}}$ an annular surface
with $\alpha^{\prime}$ non-separating and $T,\Sigma,P,Q$ surfaces of type
$(T),(\Sigma),(P),(Q)$ respectively. At the end of each of the preceding
subsections, topological properties for the link of each of the eight classes
of vertices were stated. These properties suffice to show that all eight
classes of links are pair-wise non-isomorphic as complexes.

By abuse of language, if $X$ is a separating meridian (resp. non-separating
meridian, annulus with separating boundary, annulus with non-separating
boundary, surface of type $\left.  \left(  T\right)  ,\left(  \Sigma \right)
,\left(  P\right)  ,\left(  Q\right)  \right)  $ we will be saying that
$\left[  X\right]  $ is a separating meridian vertex (resp. non-separating
meridian vertex, an annular vertex with separating boundary. an annular vertex
with non-separating boundary, a vertex of type $\left.  \left(  T\right)
,\left(  \Sigma \right)  ,\left(  P\right)  ,\left(  Q\right)  \right)  .$

Property ($sM$) characterizes the link of a vertex which is a separating
meridian. By saying "characterizes" we mean that the link of any separating
meridian has property ($sM$) and the link of any other type of vertex does not
have property ($sM$). This means that if $f\in Aut\left(  \mathcal{I}\left(
M\right)  \right)  $ and $D$ is a separating meridian then $f\left(  \left[
D\right]  \right)  $ is an isotopy class of separating meridians. Having seven
classes of links left to distinguish, property ($nsM$-$2$) characterizes the
link of a vertex which is a non-separating meridian (each of the rest $6$
classes of links contains simple closed cycles of length $3$). Thus, if $f\in
Aut\left(  \mathcal{I}\left(  M\right)  \right)  $ and $D^{\prime}$ is a
non-separating meridian, $f\left(  \left[  D^{\prime}\right]  \right)  $ is an
isotopy class of non-separating meridians. In particular $f\left(
\mathcal{D}\left(  M\right)  \right)  =\mathcal{D}\left(  M\right)  .$

Property ($nsP$) characterizes the link of a surface of type $\left(
P_{3}\right)  $ i.e., when $P$ is a pair of pants with all three boundary
components being non-separating. Combining the above mentioned properties it
can be easily seen that the same holds for all remaining types of vertices:
the link of a surface of type $\left(  P\right)  $ and the link of a surface
of type $\left(  Q\right)  $ are the only ones containing finitely many
vertices of infinite valence (see properties ($P$-$1$), ($Q$-$1$) and
properties ($sA$-$2$), ($nsA$-$2$), ($T$-$1$), ($\Sigma$-$1$)). Properties
($P$-$2$) and ($Q$-$2.1$), ($Q$-$2.2$), ($Q$-$2.3$) suffice to distinguish
between vertices of type $\left(  P\right)  $ and $\left(  Q\right)  .$ This
will be explained later.

The remaining four classes of vertices, namely, annulus with separating
boundary, annulus with non-separating boundary and surfaces of type $\left(
T\right)  ,\left(  \Sigma \right)  $ can be characterized by looking at the
valence of their vertices: property ($T$-$1$), characterizes vertices of type
$\left(  T\right)  ,$ existence of vertices of valence $7$ or $8$ (see
property ($nsA$-$1$)) characterizes annular vertices with non-separating
boundary and properties ($sA$-$3$), ($\Sigma$-$1$) suffice to distinguish
between vertices of type $\left(  \Sigma \right)  $ and annular vertices with
separating boundary.

To complete the proof of Proposition \ref{non_iso} it remains to show how
properties ($P$-$2$) and ($Q$-$2.1$), ($Q$-$2.2$), ($Q$-$2.3$) can be used in
order to distinguish between vertices of type $\left(  P\right)  $ and
$\left(  Q\right)  .$ Let $\left[  Q_{\alpha,\beta}\right]  $ be a vertex with
$Q_{\alpha,\beta}$ being a surface of type $\left(  Q\right)  $ and $\left[
P\right]  $ a vertex of type $\left(  P\right)  .$ We will consider three
cases according to whether $\left[  \Sigma_{\alpha}\right]  $ and/or $\left[
\Sigma_{\beta}\right]  $ exist in $Lk\left(  \left[  Q_{\alpha,\beta}\right]
\right)  .$

\begin{itemize}
\item If neither $\left[  \Sigma_{\alpha}\right]  $ nor $\left[  \Sigma
_{\beta}\right]  $ exist in $Lk\left(  \left[  Q_{\alpha,\beta}\right]
\right)  $ then by ($Q$-$2.1$), $Lk\left(  \left[  Q_{\alpha,\beta}\right]
\right)  $ does not contain $K_{6}$ as a subgraph whereas $Lk\left(  \left[
P\right]  \right)  $ does.

\item If exactly one of the vertices $\left[  \Sigma_{\alpha}\right]  $ and
$\left[  \Sigma_{\beta}\right]  $ exists in $Lk\left(  \left[  Q_{\alpha
,\beta}\right]  \right)  $ then, by property ($Q$-$2.2$) any two subgraphs of
$Lk\left(  \left[  Q_{\alpha,\beta}\right]  \right)  $ isomorphic to $K_{6}$
have $3$ vertices in common. This is not true for the link of a surface of
type $\left(  P\right)  :$ denote by $K_{i}$ the subgraph of $Lk\left(
\left[  P\right]  \right)  $ consisting of the vertices $\left[  S_{\alpha
}\right]  ,$ $\left[  S_{\beta}\right]  ,$ $\left[  S_{\alpha_{i}}\right]  ,$
$\left[  P_{\beta,\alpha_{i}}\right]  ,$ $\left[  \Sigma_{\alpha_{i}}\right]
,$ $\left[  Q_{\alpha,\alpha_{i}}\right]  $ (see property ($P$-$2$)). Then, by
choosing $i^{\prime}\neq i$ so that the curves $\alpha_{i}$ and $\alpha
_{i^{\prime}}$ intersect, we have two subgraphs $K_{i},$ $K_{i^{\prime}}$ of
$Lk\left(  \left[  P\right]  \right)  $ which do not have $3$ vertices in common

\item If both $\left[  \Sigma_{\alpha}\right]  $ and $\left[  \Sigma_{\beta
}\right]  $ exist in $Lk\left(  \left[  Q_{\alpha,\beta}\right]  \right)  $
then by choosing $i,$ $i^{\prime},$ $i^{\prime \prime}$ so that the isotopy
classes of the curves $\alpha_{i},$ $\alpha_{i^{\prime}},$ $\alpha
_{i^{\prime \prime}}$ are pair-wise distinct we obtain three subgraphs $K_{i},$
$K_{i^{\prime}},$ $K_{i^{\prime \prime}}$ of $Lk\left(  \left[  P\right]
\right)  $ which have the vertices $\left[  S_{\alpha}\right]  ,$ $\left[
S_{\beta}\right]  $ in common and all other vertices ($12$ of them in total)
are pair-wise distinct. This cannot be done in $Lk\left(  \left[
Q_{\alpha,\beta}\right]  \right)  $ because by property ($Q$-$2.2$) any
subgraph of $Lk\left(  \left[  Q\right]  \right)  $ isomorphic to $K_{6}$
contains the vertices $\left[  S_{\alpha}\right]  ,$ $\left[  S_{\beta
}\right]  $ and exactly one of the vertices $\left[  \Sigma_{\alpha}\right]  $
and $\left[  \Sigma_{\beta}\right]  .$
\end{itemize}

This completes the proof of Proposition \ref{non_iso}. Moreover, we have shown
the following Corollary.

\begin{corollary}
\label{non_iso_more}If $I$ is a surface of type $\left(  T\right)  $ (resp. of
type $\left(  \Sigma \right)  ,\left(  P\right)  $,$\left(  Q\right)  )$ and
$f\in Aut\left(  \mathcal{I}\left(  M\right)  \right)  $ then $f\left(
\left[  I\right]  \right)  $ is an isotopy class containing surfaces of type
$\left(  T\right)  $ (resp. of type $\left(  \Sigma \right)  ,\left(  P\right)
$,$\left(  Q\right)  ).$
\end{corollary}

We conclude this section by establishing hyperbolicity for $\mathcal{I}\left(
M\right)  .$

\begin{proposition}
If $M$ is a handlebody of genus $n\geq2,$ the complex $\mathcal{I}\left(
M\right)  $ is $\delta-$hyperbolic in the sense of Gromov.
\end{proposition}

\begin{proof}
As far as hyperbolicity is concerned, the 1-skeleton $\mathcal{I}\left(
M\right)  ^{(1)}$ of $\mathcal{I}\left(  M\right)  $ is relevant.
$\mathcal{I}\left(  M\right)  ^{(1)}$ is endowed with the combinatorial metric
so that each edge has length $1.$ Apparently, we have an embedding
\[
i:\mathcal{C}\left(  \partial M\right)  ^{(1)}\hookrightarrow \mathcal{I}%
\left(  M\right)  ^{(1)}%
\]
with $i:\mathcal{C}\left(  \partial M\right)  ^{(1)}=\mathcal{D}\left(
M\right)  ^{(1)}\cup \mathcal{A}\left(  M\right)  ^{(1)}$ where the superscript
$^{(1)}$ always denotes $1-$skeleton. We claim that this embedding is
isometric. Indeed, if $\left[  \alpha_{1}\right]  ,\left[  \alpha_{2}\right]
$ are distinct vertices with distance $d_{\mathcal{C}}\left(  \left[
\alpha_{1}\right]  ,\left[  \alpha_{2}\right]  \right)  $ in $\mathcal{C}%
\left(  \partial M\right)  ^{(1)}$ then the distance $d_{\mathcal{I}}\left(
i\left(  \left[  \alpha_{1}\right]  \right)  ,i\left(  \left[  \alpha
_{2}\right]  \right)  \right)  $ cannot be smaller. For, if $\left[
S_{0}\right]  =i\left(  \left[  \alpha_{1}\right]  \right)  ,\left[
S_{1}\right]  ,\ldots,\left[  S_{k}\right]  =i\left(  \left[  \alpha
_{2}\right]  \right)  $ is a sequence of vertices which gives rise to a
geodesic in $\mathcal{I}\left(  M\right)  ^{(1)}$ of length less than
$d_{\mathcal{C}}\left(  \left[  \alpha_{1}\right]  ,\left[  \alpha_{2}\right]
\right)  ,$ equivalently,
\[
d_{\mathcal{I}}\left(  i\left(  \left[  \alpha_{1}\right]  \right)  ,i\left(
\left[  \alpha_{2}\right]  \right)  \right)  =k<d_{\mathcal{C}}\left(  \left[
\alpha_{1}\right]  ,\left[  \alpha_{2}\right]  \right)
\]
then for each $j=1,2,\ldots,k-1$ consider $\partial S_{j}$ to be any boundary
component of $S_{j}.$ It is clear that $\partial S_{j}$ is disjoint from
$\partial S_{j-1}$ and $\partial S_{j+1}.$ Therefore, the sequence $\left[
\alpha_{1}\right]  ,\left[  \partial S_{1}\right]  ,\ldots,\left[  \partial
S_{k-1}\right]  ,\left[  \alpha_{2}\right]  $ is a segment of length $k$ with
$k<d_{\mathcal{C}}\left(  \left[  \alpha_{1}\right]  ,\left[  \alpha
_{2}\right]  \right)  $, a contradiction.

For any vertex $\left[  \Sigma \right]  $ in $\mathcal{I}\left(  M\right)
^{(1)}\setminus \mathcal{D}\left(  M\right)  ^{(1)}\cup \mathcal{A}\left(
M\right)  ^{(1)}$ we may find an annular vertex$,$ namely, $S_{\partial \Sigma
}$ where $\partial \Sigma$ is any component of the boundary of $\Sigma,$ which
is connected by an edge with $\left[  \Sigma \right]  .$ Thus, $\mathcal{I}%
\left(  M\right)  ^{(1)}$ is within bounded distance from $i\left(
\mathcal{C}\left(  \partial M\right)  ^{(1)}\right)  .$ Since $\mathcal{C}%
\left(  \partial M\right)  ^{(1)}$ is $\delta-$hyperbolic in the sense of
Gromov, so is $\mathcal{I}\left(  M\right)  ^{(1)}.$
\end{proof}

\section{Proof of the Main Theorem\label{main}}

\bigskip Let $M$ be a handlebody of genus $n=2.$ If $F$ is a
self-homeomorphism of $M,$ it is clear that $F$ sends incompressible surfaces
to incompressible surfaces, isotopic surfaces to isotopic surfaces and, hence,
isotopy classes of incompressible surfaces to isotopy classes of
incompressible surfaces. In other words $F$ induces a morphism denoted by
$A\left(  F\right)  $ of the complex $Aut\left(  \mathcal{I}\left(  M\right)
\right)  $ given by
\[
A\left(  F\right)  \left[  S\right]  :=\left[  F\left(  S\right)  \right]  .
\]
As $F$ is invertible this morphism is an automorphism. Finally, if $F$ is
isotopic to $F^{\prime}$ then $A\left(  F\right)  =A\left(  F^{\prime}\right)
$ since $\mathcal{I}\left(  M\right)  $ is a flag complex defined up to
isotopy. Therefore we have a well defined map
\[
A:\mathcal{MCG}\left(  M\right)  \rightarrow Aut\left(  \mathcal{I}\left(
M\right)  \right)
\]
where $\mathcal{MCG}\left(  M\right)  $ denotes the group of isotopy classes
of self-homeomorphisms of $M.$ In an identical way the map
\[
A_{0}:\mathcal{MCG}\left(  M\right)  \rightarrow Aut\left(  \mathcal{I}%
_{0}\left(  M\right)  \right)
\]
is well defined, where $\mathcal{I}_{0}\left(  M\right)  $ the subcomplex of
$\mathcal{I}\left(  M\right)  $ consisting of all vertices of infinite valence
in $\mathcal{I}\left(  M\right)  ,$ i.e.
\[
\mathcal{I}_{0}\left(  M\right)  =\mathcal{I}\left(  M\right)  \setminus
\left \{  \left[  S\right]  \bigm \vert S\mathrm{\ is\ of\ type\ }\left(
P_{3}\right)  \right \}
\]

\begin{theorem}
\label{main_theorem} The map $A_{0}:\mathcal{MCG}\left(  M\right)  \rightarrow
Aut\left(  \mathcal{I}_{0}\left(  M\right)  \right)  $ is onto and has a
$\mathbb{Z}_{2}-$kernel. The map $A:\mathcal{MCG}\left(  M\right)  \rightarrow
Aut\left(  \mathcal{I}\left(  M\right)  \right)  $ is injective and
$Aut\left(  \mathcal{I}\left(  M\right)  \right)  $ contains non-geometric elements.
\end{theorem}

We will need the following lemma.

\begin{lemma}
\label{l1}Let $M=H_{2}$ be the handlebody of genus $2.$ If $f$ $\in Aut\left(
\mathcal{I}\left(  M\right)  \right)  $ and $f|_{\mathcal{C}\left(  \partial
M\right)  }=id_{\mathcal{C}\left(  \partial M\right)  }$ then $f\left(
\left[  S\right]  \right)  =\left[  S\right]  $ for any vertex $\left[
S\right]  \in \mathcal{I}_{0}\left(  M\right)  .$ If $\left[  S\right]
=\left[  P_{\alpha,\beta,\gamma}^{+}\right]  \in \mathcal{I}\left(  M\right)
\setminus \mathcal{I}_{0}\left(  M\right)  $ then either, $f\left(  \left[
S\right]  \right)  =\left[  S\right]  $ or, $f\left(  \left[  S\right]
\right)  =\left[  P_{\alpha,\beta,\gamma}^{-}\right]  .$
\end{lemma}

\begin{proof}
We first show that $f$ fixes all vertices of type $\left(  \Sigma \right)  .$
Recall that by Lemma \ref{non_iso_more}, the image of a vertex of type
$\left(  \Sigma \right)  $ is a vertex of type $\left(  \Sigma \right)  .$ Let
$\alpha,\beta$ be two non-separating non-isotopic (i.e. $\left[
\alpha \right]  \neq \left[  \beta \right]  $) curves and $\Sigma_{\alpha}%
,\Sigma_{\beta}$ the corresponding surfaces of type $(\Sigma)$. If $\left[
a\right]  \cap \left[  \beta \right]  \neq \varnothing$ then the annular surface
$\left[  S_{\alpha}\right]  $ is a vertex for which the edge $\left(  \left[
S_{\alpha}\right]  ,\left[  \Sigma_{\alpha}\right]  \right)  $ exists whereas
the edge $\left(  \left[  S_{\alpha}\right]  ,\left[  \Sigma_{\beta}\right]
\right)  $ does not. As $f$ is assumed to fix all annular vertices $f$ cannot
map, in this case, $\left[  \Sigma_{\alpha}\right]  $ onto $\left[
\Sigma_{\beta}\right]  .$ If $\left[  a\right]  \cap \left[  \beta \right]
=\varnothing$ then choose a curve $\gamma$ such that $\left[  \alpha \right]
\cap \left[  \gamma \right]  =\varnothing$ and $\left[  \beta \right]
\cap \left[  \gamma \right]  \neq \varnothing$ (e.g. $\gamma$ may be $\beta^{2}%
$). Then $\left[  S_{\gamma}\right]  $ is an annular vertex for which the edge
$\left(  \left[  S_{\gamma}\right]  ,\left[  \Sigma_{\alpha}\right]  \right)
$ exists whereas the edge $\left(  \left[  S_{\gamma}\right]  ,\left[
\Sigma_{\beta}\right]  \right)  $ does not. Thus, for any $\alpha,\beta$ with
$\left[  \alpha \right]  \neq \left[  \beta \right]  $, $f$ cannot map $\left[
\Sigma_{\alpha}\right]  $ onto $\left[  \Sigma_{\beta}\right]  .$ Thus $f$
must fix all vertices of type $\left(  \Sigma \right)  .$

We proceed to show that $f$ must fix all vertices of type $\left(  T\right)
.$ Again by Lemma \ref{non_iso_more} the image of a vertex of type $\left(
T\right)  $ is a vertex of type $\left(  T\right)  .$ Let $T_{\alpha,+}$ be a
surface of type $\left(  T\right)  $ for an arbitrary separating curve
$\alpha$ (the proof for $T_{\alpha,-}$ will be identical). Let $f\left(
\left[  T_{\alpha,+}\right]  \right)  =\left[  T_{\beta,+}\right]  $ for some
separating curve $\beta$ non-isotopic to $\alpha.$ The non-empty intersection
$\left[  a\right]  \cap \left[  \beta \right]  \neq \varnothing$ implies that
$\beta$ contains subarcs of the form $\beta_{1},\beta_{2}$ shown in Figure
\ref{fig11}. Choose a curve $\alpha_{i_{0}}$ intersecting $\beta_{1}$ or
$\beta_{2}$ and not intersecting $\alpha.$ Then $\left[  S_{\alpha_{i_{0}}%
}\right]  \in Lk\left(  \left[  T_{\alpha,+}\right]  \right)  $ and the edge
$\left(  \left[  T_{\alpha,+}\right]  ,\left[  S_{\alpha_{i_{0}}}\right]
\right)  $ exists whereas $\left(  \left[  T_{\beta,+}\right]  ,\left[
S_{\alpha_{i_{0}}}\right]  \right)  $ does not. As $f$ is assumed to fix all
annular vertices, it cannot map $\left[  T_{\alpha,+}\right]  $ onto $\left[
T_{\beta,+}\right]  .$ Similarly for $\left[  T_{\beta,-}\right]  .$ It
remains to verify that $f$ cannot map $\left[  T_{\alpha,+}\right]  $ onto
$\left[  T_{\alpha,-}\right]  .$ In fact, this is not possible as for
arbitrary index $j_{0},$ the vertex $\left[  \Sigma_{\alpha_{j_{0}}^{-}%
}\right]  $ belongs to $Lk\left(  \left[  T_{\alpha,+}\right]  \right)  ,$
does not belong to $Lk\left(  \left[  T_{\alpha,-}\right]  \right)  $ and the
edge $\left(  \left[  T_{\alpha,+}\right]  ,\left[  \Sigma_{\alpha_{j_{0}}%
^{-}}\right]  \right)  $ exists whereas $\left(  \left[  T_{\alpha,-}\right]
,\left[  \Sigma_{\alpha_{j_{0}}^{-}}\right]  \right)  $ does not.

We next examine vertices of type $\left(  P_{3}\right)  $ i.e., a pair of
pants $P_{\alpha_{1},\alpha_{2},\alpha_{3}}^{+}$ with $\alpha_{1},\alpha
_{2},\alpha_{3}$ non-separating boundary curves. Since the link $Lk\left(
\left[  P_{\alpha_{1},\alpha_{2},\alpha_{3}}^{+}\right]  \right)  $ contains
finitely many vertices, $f\left(  \left[  P_{\alpha_{1},\alpha_{2},\alpha_{3}%
}^{+}\right]  \right)  $ must necessarily be a vertex $\left[  P_{\beta
_{1},\beta_{2},\beta_{3}}\right]  $ for a pair of pants $P_{\beta_{1}%
,\beta_{2},\beta_{3}}$ with~$\beta_{1},\beta_{2},\beta_{3}$ non-separating
curves in $\partial H_{2}.$ If $\left[  \alpha_{i}\right]  \cap \left[
\beta_{j}\right]  \neq \varnothing$ for some $i,j\in \left \{  1,2,3\right \}  $
then the edge $\left(  \left[  S_{\alpha_{i}}\right]  ,\left[  P_{\alpha
_{1},\alpha_{2},\alpha_{3}}^{+}\right]  \right)  $ exists whereas the edge
$\left(  \left[  S_{\alpha_{i}}\right]  ,\left[  P_{\beta_{1},\beta_{2}%
,\beta_{3}}\right]  \right)  $ does not. It follows that, up to a change of
enumeration, $\alpha_{i}$ is isotopic to $\beta_{i}$ for $i=1,2,3.$ Therefore,
$P_{\beta_{1},\beta_{2},\beta_{3}}$ is isotopic to, either $P_{\alpha
_{1},\alpha_{2},\alpha_{3}}^{+}$ or, $P_{\alpha_{1},\alpha_{2},\alpha_{3}}%
^{-}.$ In other words, either $f\left(  \left[  P_{\alpha_{1},\alpha
_{2},\alpha_{3}}^{+}\right]  \right)  =\left[  P_{\alpha_{1},\alpha_{2}%
,\alpha_{3}}^{+}\right]  $ or, $f\left(  \left[  P_{\alpha_{1},\alpha
_{2},\alpha_{3}}^{+}\right]  \right)  =\left[  P_{\alpha_{1},\alpha_{2}%
,\alpha_{3}}^{-}\right]  .$ \newline Let now $\left[  P_{\alpha,\beta}\right]
$ be a vertex of type $\left(  P\right)  $ (with $\alpha$ non-separating and
$\beta$ separating). If $\left[  P_{\gamma,\beta^{\prime}}\right]  $ is a
surface of type $\left(  P\right)  $ with $\left[  \beta^{\prime}\right]
\cap \left[  \beta \right]  \neq \varnothing$ (and $\gamma$ arbitrary) then the
edge $\left(  \left[  T_{\beta,\pm}\right]  ,\left[  P_{\alpha,\beta}\right]
\right)  $ exists whereas the edge $\left(  \left[  T_{\beta,\pm}\right]
,\left[  P_{\gamma,\beta^{\prime}}\right]  \right)  $ does not. As $f$ fixes
$\left[  T_{\beta,+}\right]  $ (or, $\left[  T_{\beta,-}\right]  $), this
shows that $f\left(  \left[  P_{\alpha,\beta}\right]  \right)  \neq \left[
P_{\gamma,\beta^{\prime}}\right]  $ for all choices of $\gamma$ provided that
$\left[  \beta^{\prime}\right]  \cap \left[  \beta \right]  \neq \varnothing.$
Similarly, using the annular vertex $\left[  S_{\alpha}\right]  $ which is
fixed by $f$ it can be seen that $f\left(  \left[  P_{\alpha,\beta}\right]
\right)  \neq \left[  P_{\alpha^{\prime},\beta}\right]  $ for any
$\alpha^{\prime}$ with $\left[  \alpha^{\prime}\right]  \cap \left[
\alpha \right]  \neq \varnothing.$ It remains to examine whether $f$ can map
$\left[  P_{\alpha,\beta}\right]  $ to a vertex $\left[  P_{\alpha^{\prime
},\beta}\right]  $ for some curve $\alpha^{\prime}$ non-isotopic to $\alpha$
with the property $\left[  \alpha^{\prime}\right]  \cap \left[  \alpha \right]
=\varnothing.$ Under these assumptions for $\alpha,\alpha^{\prime}$ it follows
that $\alpha,\alpha^{\prime}$ belong to distinct components of $\partial
M\setminus \beta.$ Choose a curve $\alpha_{i_{0}}$ belonging to the component
of $\partial M\setminus \beta$ containing $\alpha^{\prime}.$ Then the edge
$\left(  \left[  S_{\alpha_{i_{0}}}\right]  ,\left[  P_{\alpha,\beta}\right]
\right)  $ exists whereas the edge $\left(  \left[  S_{\alpha_{i_{0}}}\right]
,\left[  P_{\alpha^{\prime},\beta}\right]  \right)  $ does not. Thus $f$ fixes
all vertices of type $\left(  P\right)  $ $.$

\begin{figure}[ptb]
\begin{center}
\includegraphics[scale=0.6]
{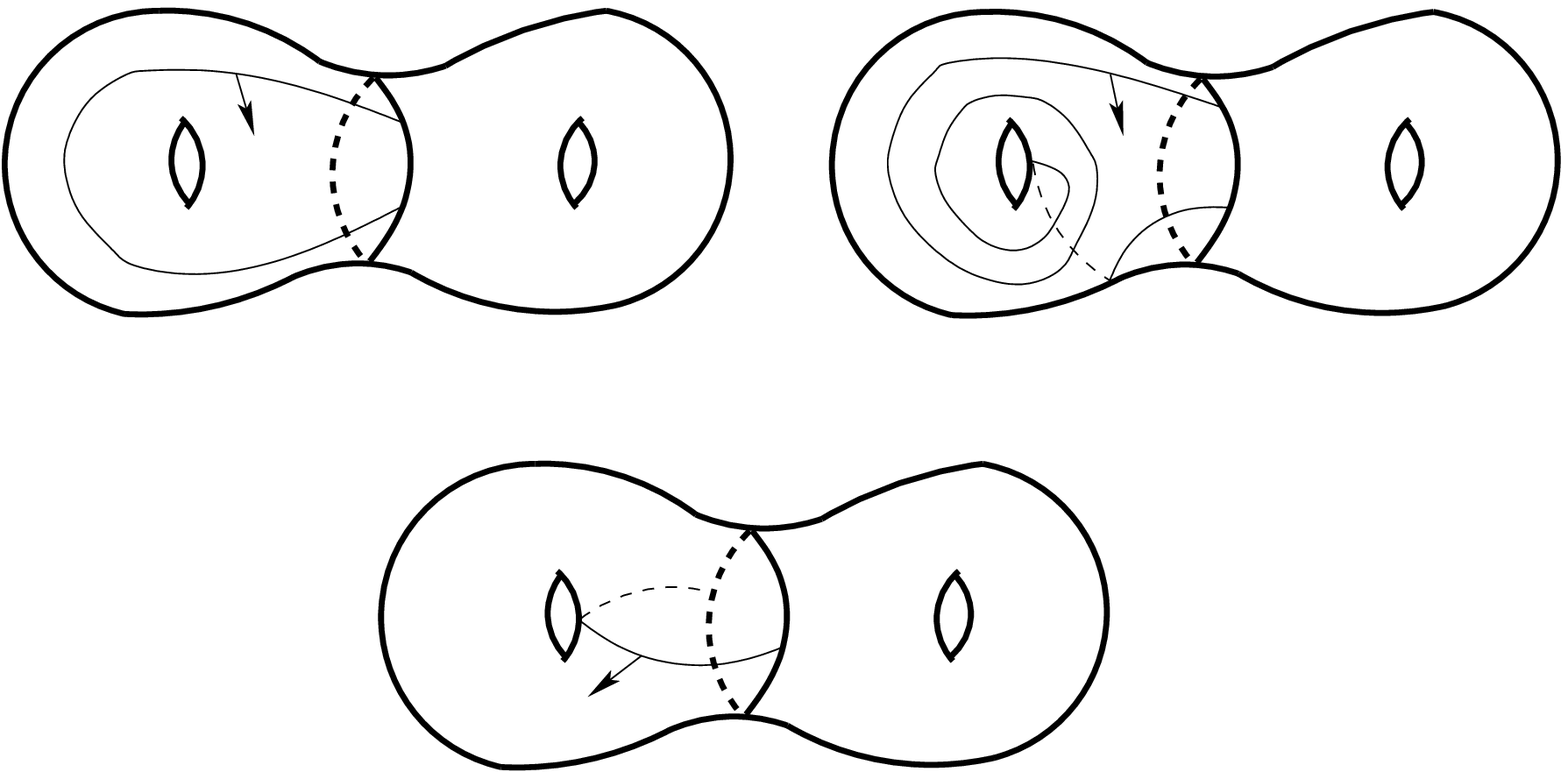}
\end{center}
\par
\begin{picture}(22,12)
\put(100,146){$\beta_1$}
\put(120,115){$\alpha$}
\put(284,116){$\alpha$}
\put(194,27){$\alpha$}
\put(270,144){$\beta_1$}
\put(159,38){$\beta_2$}
\end{picture}
\caption{$\,$}%
\label{fig11}%
\end{figure}

Finally, let $Q_{\alpha,\beta},Q_{\gamma,\delta}$ be two non-isotopic (i.e.
$\left[  Q_{\alpha,\beta}\right]  \neq \left[  Q_{\gamma,\delta}\right]  $ )
surfaces with $f\left(  \left[  Q_{\alpha,\beta}\right]  \right)  =\left[
Q_{\gamma,\delta}\right]  .$ Apparently, the edges $\left(  \left[  S_{\alpha
}\right]  ,\left[  Q_{\alpha,\beta}\right]  \right)  $, $\left(  \left[
S_{\beta}\right]  ,\left[  Q_{\alpha,\beta}\right]  \right)  $ exist in
$Lk\left(  \left[  Q_{\alpha,\beta}\right]  \right)  $ and as $f$ is assumed
to fix all annular vertices, it follows that $\left(  \left[  S_{\alpha
}\right]  ,\left[  Q_{\gamma,\delta}\right]  \right)  $, $\left(  \left[
S_{\beta}\right]  ,\left[  Q_{\gamma,\delta}\right]  \right)  $ exist in
$Lk\left(  \left[  Q_{\gamma,\delta}\right]  \right)  .$ Thus, $\gamma,\delta$
are two non-separating curves and each of them does not intersect both
$\alpha$ and $\beta.$ In other words, $\gamma,\delta$ are two disjoint curves
in the sphere $Q_{\alpha,\beta}$ with four holes. Assume one of them, say
$\gamma,$ is not isotopic to neither $\alpha$ nor $\beta.$ Consider the
(essential) separating curve $\overline{\alpha}$ obtained by joining two
copies of $\alpha$ with a simple arc which intersects $\gamma$ but not
$\beta.$ Then, for the vertex $\left[  P_{\overline{\alpha},\alpha}\right]  $
we have that the edge $\left(  \left[  Q_{\alpha,\beta}\right]  ,\left[
P_{\overline{\alpha},\alpha}\right]  \right)  $ exists whereas the edge
$\left(  \left[  Q_{\gamma,\delta}\right]  ,\left[  P_{\overline{\alpha}%
,a}\right]  \right)  $ does not exist because $\overline{\alpha}$ intersects
$\gamma.$ This is not possible since we assumed that $f\left(  \left[
Q_{\alpha,\beta}\right]  \right)  =\left[  Q_{\gamma,\delta}\right]  $ and
$\left[  P_{\overline{\alpha},a}\right]  $ is fixed by $f,$ as shown above.
Therefore, we may assume that $\gamma$ is isotopic to $\alpha$ and similarly
we obtain that $\delta$ is isotopic to $\beta.$ It follows that $\left[
Q_{\alpha,\beta}\right]  =\left[  Q_{\gamma,\delta}\right]  ,$ a
contradiction. Hence, $f$ fixes all vertices of type $\left(  Q\right)  $ and
this completes the proof of the Lemma.
\end{proof}\\[1mm]

\begin{proof}
[Proof of Theorem \ref{main_theorem}]We will use the corresponding result for
surfaces, see \cite{I3},\cite{Luo}, which applies to the boundary of the handlebody
$M=H_{2}.$ The map $\mathcal{MCG}\left(  \partial M\right)  \rightarrow
Aut\left(  \mathcal{C}\left(  \partial M\right)  \right)  $ will be denoted
again by $A_{0}.$ By abuse of language, using the identification
\[
\mathcal{M}\left(  M\right)  \cup \mathcal{A}\left(  M\right)
\longleftrightarrow \mathcal{C}\left(  \partial M\right)  ,
\]
we will be viewing the complex $\mathcal{C}\left(  \partial M\right)  $ as a
subcomplex of $\mathcal{I}_{0}\left(  M\right)  .$\newline

We first show that every $f\in Aut\left(  \mathcal{I}_{0}\left(  M\right)
\right)  $ is geometric. By Proposition \ref{non_iso} we know that $f\left(
\mathcal{A}\left(  M\right)  \right)  =\mathcal{A}\left(  M\right)  $ and
$f\left(  \mathcal{M}\left(  M\right)  \right)  =\mathcal{M}\left(  M\right)
.$ In particular, $f\left(  \mathcal{C}\left(  \partial M\right)  \right)
=\mathcal{C}\left(  \partial M\right)  .$ The restriction $f|_{\mathcal{C}%
\left(  \partial M\right)  }$ of $f$ on $\mathcal{C}\left(  \partial M\right)
$ induces an automorphism of $\mathcal{C}\left(  \partial M\right)  $ which by
the analogous result for surfaces (see \cite{I3},\cite{Luo}) is geometric, 
that is, there exists a homeomorphism
\[
F_{\partial M}:\partial M\rightarrow \partial M
\]
such that $A_{0}\left(  F_{\partial M}\right)  =f|_{\mathcal{C}\left(
\partial M\right)  }.$ As $f|_{\mathcal{C}\left(  \partial M\right)  }$ maps
$\mathcal{M}\left(  M\right)  $ to $\mathcal{M}\left(  M\right)  ,$
$F_{\partial M}$ sends meridian boundaries to meridian boundaries. It follows
that $F_{\partial M}$ extends to a homeomorphism $F:M\rightarrow M.$ We know
that $A_{0}\left(  F\right)  =f$ on $\mathcal{C}\left(  \partial M\right)  $
and we must show that $A_{0}\left(  F\right)  =f$ on $\mathcal{I}\left(
M\right)  .$ This follows from Lemma \ref{l1} which completes the proof every
$f\in Aut\left(  \mathcal{I}_{0}\left(  M\right)  \right)  $ is geometric.
\newline Let $f\in Aut\left(  \mathcal{I}_{0}\left(  M\right)  \right)  .$
Since $A_{0}$ is shown to be onto, there exists a homeomorphism
$F:M\rightarrow M$ such that $A\left(  \left[  F\right]  \right)  =f.$ This
implies that $f\left(  \mathcal{M}\left(  M\right)  \right)  =\mathcal{M}%
\left(  M\right)  $ and $f\left(  \mathcal{A}\left(  M\right)  \right)
=\mathcal{A}\left(  M\right)  .$ In particular, $f$ restricted to
$\mathcal{C}\left(  \partial M\right)  \equiv \mathcal{M}\left(  M\right)
\cup \mathcal{A}\left(  M\right)  $ induces an automorphism $\overline{f}$ of
the complex of curves $\mathcal{C}\left(  \partial M\right)  .$ By \cite{I3}, \cite{Luo},
there exists a homeomorphism $F_{\partial M}:\partial M\rightarrow \partial M$
such that $A_{0} \left(  F_{\partial M}\right)  =\overline{f}.$ Such a
homeomorphism is not unique because the map
\[
\mathcal{MCG}\left(  \partial M\right)  \rightarrow Aut\left(  \mathcal{C}%
\left(  \partial M\right)  \right)
\]
has a $\mathbb{Z}_{2}-$kernel generated by an involution of $\partial M.$
However, any homeomorphism of $\partial M$ which extends to $M$ it does so
uniquely (see, for example, \cite[Theorem 3.7 p.94]{F-M}), and therefore the
map
\[
\mathcal{MCG}\left(  M\right)  \rightarrow Aut\left(  \mathcal{I}_{0}\left(
M\right)  \right)
\]
also has a $\mathbb{Z}_{2}-$kernel.

We now show that the map $\mathcal{MCG}\left(  M\right)  \rightarrow
Aut\left(  \mathcal{I}\left(  M\right)  \right)  $ is injective but not
surjective. First observe that if $I:M\rightarrow M$ is an involution then
$A(I)$ fixes $\mathcal{C}\left(  \partial M\right)  $ and, by Lemma \ref{l1},
$A(I)$ fixes $\mathcal{I}_{0}\left(  M\right)  .$ Moreover, it fixes all
elements $\left[  S\right]  =\left[  P_{\alpha,\beta,\gamma}^{+}\right]  $ in
$\mathcal{I}\left(  M\right)  \setminus \mathcal{I}_{0}\left(  M\right)  $ for
which $P_{\alpha,\beta,\gamma}^{+},$ $P_{\alpha,\beta,\gamma}^{+}$ are
isotopic. If $\left[  S\right]  =\left[  P_{\alpha,\beta,\gamma}^{+}\right]  $
is any vertex of type $\left(  P_{3}\right)  $ such that the pairs of pants
$P_{\alpha,\beta,\gamma}^{+},$ $P_{\alpha,\beta,\gamma}^{-}$ are non isotopic
(see Remark \ref{r5}) then $I\left(  P_{\alpha,\beta,\gamma}^{+}\right)  $ is
isotopic to $P_{\alpha,\beta,\gamma}^{-}$ which means that $A(I)$ interchanges
$\left[  P_{\alpha,\beta,\gamma}^{+}\right]  ,$ $\left[  P_{\alpha
,\beta,\gamma}^{-}\right]  $ for any such vertex $\left[  S\right]  =\left[
P_{\alpha,\beta,\gamma}^{+}\right]  .$ In particular, this shows that $A$ is injective.

To define a non-geometric element $g\in Aut\left(  \mathcal{I}\left(
M\right)  \right)  $ pick non-separating curves $\alpha_{0},$ $\beta_{0},$
$\gamma_{0}$ in $\partial M$ such that $P_{\alpha_{0},\beta_{0},\gamma_{0}%
}^{+},$ $P_{\alpha_{0},\beta_{0},\gamma_{0}}^{-}$ are non-isotopic and set

\begin{itemize}
\item $g\left(  \left[  P_{\alpha_{0},\beta_{0},\gamma_{0}}^{+}\right]
\right)  =\left[  P_{\alpha_{0},\beta_{0},\gamma_{0}}^{-}\right]  ,$ $g\left(
\left[  P_{\alpha_{0},\beta_{0},\gamma_{0}}^{-}\right]  \right)  =\left[
P_{\alpha_{0},\beta_{0},\gamma_{0}}^{+}\right]  $ and

\item $g\left(  \left[  S\right]  \right)  =\left[  S\right]  $ for all
$\left[  S\right]  \in \mathcal{I}\left(  M\right)  \setminus \left \{  \left[
P_{\alpha_{0},\beta_{0},\gamma_{0}}^{+}\right]  ,\left[  P_{\alpha_{0}%
,\beta_{0},\gamma_{0}}^{-}\right]  \right \}  .$
\end{itemize}

Assuming that there exists a homeomorphism $G:M\rightarrow M$ such that
$A\left(  \left[  G\right]  \right)  =g,$ then, since $g$ fixes $\mathcal{I}%
_{0}\left(  M\right)  ,$ it follows that $G$ is in the kernel of $A_{0}.$
Thus, $G$ is either, the identity or, an involution. Apparently, $G$ cannot be
the identity as $g$ is non trivial and $G$ cannot be an involution because the
image under $A$ of an involution interchanges all pairs of vertices $\left[
P_{\alpha,\beta,\gamma}^{+}\right]  ,\left[  P_{\alpha,\beta,\gamma}%
^{-}\right]  $ for which $P_{\alpha,\beta,\gamma}^{+},$ $P_{\alpha
,\beta,\gamma}^{-}$ are non isotopic.
\end{proof}

\bigskip \textit{Acknowledgment.} We would like to thank A. Papadopoulos for
asking us whether a complex, analogous to the complex of curves, can be
defined in order to encode the homeomorphisms of a handlebody by the
automorphism group of the complex.

\end{document}